\documentclass[12pt, leqno]{article}
\usepackage{tikz}
\usepackage{float}
\usepackage{amsfonts}
\usepackage{amsmath}
\usepackage{amssymb}
\usepackage{amscd}
\usepackage{mathrsfs}
\usepackage{color}

\begin{document}
\renewcommand{\thefootnote}{\fnsymbol{footnote}}
\begin{center}
{\bf \Large On $3$-dimensional foliated dynamical systems and Hilbert type reciprocity law}
\end{center} 

\vspace{0.2cm}

\begin{center}
Junhyeong Kim, Masanori Morishita, Takeo Noda and Yuji Terashima
\end{center}

\begin{center}
{\em Dedicated to Professor Christopher Deninger}
\end{center}

\vspace{0.2cm}


{\small
Abstract:  We show some fundamental results concerning $3$-dimensional foliated dynamical systems (FDS$^3$ for short) introduced by Deninger. Firstly, we give a decomposition theorem for an FDS$^3$, which yields a classification of  FDS$^3$'s. Secondly,  for each type of the classification, we construct concrete examples of FDS$^3$'s.  Finally, by using the integration theory for smooth Deligne cohomology, we introduce geometric analogues of local symbols and show a Hilbert type reciprocity law for an FDS$^3$. Our results answer the question posed by Deninger.}

\vspace{1.2cm}

\begin{center}
{\bf Introduction}
\end{center}

In his monumental works (cf. [D1] $\sim$ [D7]), Deninger considered the notion of a foliated dynamical system, namely, a smooth manifold equipped with $1$-codimensional foliation and transverse flow satisfying certain conditions, as a geometric analogue of an arithmetic scheme (see also [Ko1], [Ko2]). In particular, $3$-dimensional foliated dynamical systems  may be regarded as analogues of arithmetic curves, where closed orbits (knots) correspond to closed points (finite primes). The purpose of this paper is to study some fundamental questions concerning $3$-dimensional foliated dynamical systems. Let us describe our main results in the following.\\

Let  $ \mathfrak{S} = (M,{\cal F}, \phi)$ be a $3$-dimensional foliated dynamical system, called an FDS$^3$ for short. Namely, $M$ is a connected, closed smooth $3$-manifold, ${\cal F}$ is a complex foliation by Riemann surfaces on $M$, $\phi$ is a smooth dynamical system on $M$. These data must satisfy the conditions: there is a set of  finitely many compact leaves ${\cal P}_{\mathfrak{S}}^{\infty} = \{L_1^{\infty}, \dots , L_r^{\infty}\}$, which may be empty, such that for any $i$ and $t$ we have $\phi^t(L_i^{\infty}) = L_i^{\infty}$ and that any orbit of $\phi$ is transverse to leaves in $M_0 := M \setminus \cup_{i=1}^r L_i^{\infty}$, and $\phi^t$ maps any leaf to a leaf for each $t$ (cf. Definition 1.5 below). Let ${\cal P}_{\mathfrak{S}}$ be the set of closed orbits in $M_0$ and let $\overline{{\cal P}_{\mathfrak{S}}} = {\cal P}_{\mathfrak{S}} \cup {\cal P}_{\mathfrak{S}}^{\infty}$. Note that ${\cal P}_{\mathfrak{S}}^{\infty}$ (resp. ${\cal P}_{\mathfrak{S}}$) may be regarded as an analogue of the set of infinite (resp. finite) primes of a global field.
\\
@1) Our first result concerns the structure of an FDS$^3$. We show a decomposition theorem for an FDS$^3$, which asserts that if $M$ is cut along non-transverse compact leaves $L_i^{\infty}$'s, each connected component enjoys a well-described structure of a foliated $3$-manifold (cf. Theorem 2.2.2 below). This yields a classification of FDS$^3$'s according to the sets of transverse and non-transverse compact leaves (cf. Corollary 2.2.4). Our decomposition theorem may remind us of the JSJ decomposition of a $3$-manifold ([JS], [Jo]). \\
@2)  Our second result is to construct concrete examples of FDS$^3$'s. Before our work, all known examples of FDS$^3$ in the literature were only a closed surface bundle over $S^1$ with the bundle foliation and the suspension flow and Alvarez-Lopez's example ([D7; page 10]). In this paper, we construct concrete examples of FDS$^3$'s for each type of our classification. Furthermore, we show that any closed smooth $3$-manifold admits a structure of an FDS$^3$, by using an open book decomposition, which is a peculiar property for the $3$-dimensional case (cf. Example 3.III-1.2 below). In view of the analogy with an arithmetic curve which has countably infinitely many finite primes, we give examples of FDS$^3$'s such that ${\cal P}_{\mathfrak{S}}$ is a countably infinite set. \\
@3) Our third result concerns the question, which was posed by Deninger,  on finding out analogies for an FDS$^3$ of the Hilbert symbol and the reciprocity law in class field theory  for a global field ([AT; Chapter 12]).   For FDS$^3$-meromorphic functions $f$ and $g$, we  introduce the local symbol $\langle f, g \rangle_{\gamma}$ along  $\gamma \in \overline{P_{\mathfrak{S}}}$ and show the reciprocity law
 $$ \sum_{\gamma \in \overline{{\cal P}_{\mathfrak{S}}} } \langle f, g \rangle_{\gamma} = 0 \;\; \mbox{mod} \; (2\pi \sqrt{-1})^3\Lambda_{\mathfrak{S}}, $$
 where $\Lambda_{\mathfrak{S}}$ is the period group of $\mathfrak{S}$ (cf. Theorem 5.2.1 below). For this, we employ the theory of smooth Deligne cohomology groups and the integration theory of Deligne cohomology classes, which were studied by Brylinski ([Br], [BM]) and by Gawedzki ([Ga]) and Gomi-Terashima ([GT], [Te]).  For $\gamma \in {\cal P}_{\mathfrak{S}}$, we give an explicit integral formula for the local symbol $\langle f, g \rangle_{\gamma}$ (cf. Theorem 5.1.2). Our method may be regarded as a generalization to an FDS$^3$ of Deligne, Bloch and Beilinson's interpretation of the tame symbol on a Riemann surface using the holomorphic Deligne cohomology ([Be], [Bl], [Dl]). Our result is also a generalization of Stelzig's work on the case of closed surface bundles over $S^1$ ([St]). \\

We note that our third result may indicate the possibility to develop an id\`{e}lic theory for FDS$^3$'s. For this line of study, we mention the recent works by Niibo and Ueki ([NU]) and by Mihara ([Mi]) on id\`{e}lic class field theory for $3$-manifolds in arithmetic topology. In fact, our original desire was to refine and deepen the analogies between knots and primes, 3-manifolds and number rings in arithmetic topology ([Mo]), in the context of FDS$^3$'s.\\

Here are the contents of this paper. In the section 1, we introduce a $3$-dimensional foliated dynamical system (FDS$^3$) and some basic notions. In Section 2,  we show a decomposition theorem for an FDS$^3$, which yields a classification of FDS$^3$'s. In the section 3,  for each type of the classification, we construct concrete examples of FDS$^3$'s. In the section 4, we  recall the theory of smooth Deligne cohomology for an FDS$^3$ and the integration theory of Deligne cohomology classes. In the section 5, we introduce  local symbols of FDS$^3$-meromorphic functions along a closed orbit or along a non-transverse compact leaf, and show a Hilbert type reciprocity law.\\

\begin{center}
{\bf 1. $3$-dimensional foliated dynamical systems}
\end{center}

In this section, following Deninger and Kopei ([D1]$\sim$[D7], [Ko1],[Ko2]), we introduce the notion of a $3$-dimensional foliated dynamical system, called an FDS$^3$ for short, and prepare some basic notions and properties.  Although  a foliated dynamical system can be introduced in any odd dimension, we consider only the $3$-dimensional case in this paper, since we are concerned with analogies with arithmetic curves. For general materials in foliation theory, we refer to [CC]  and [Ta].  \\

 We begin to recall the notion of a complex foliation by Riemann surfaces (cf. [Gh]).\\
\\
{\bf Definition 1.1.} Let $M$ be a smooth $3$-manifold whose boundary $\partial M$ may be non-empty.  A {\em $2$-dimentional smooth foliation} ${\cal F}$ on $M$  is defined by a  family of immersed $2$-dimensional manifolds $\{ L_a\}_{a \in A}$ satisfying the following conditions:\\
(1) $ M = \sqcup_{a \in A} L_a $ (disjoint union).  \\
(2) there is a system of foliated local coordinates $(U_i; \varphi_i)_{i \in I}$, where $U_i$ is a foliated open subset of $M$ and $\varphi_i : U_i \stackrel{\approx}{\rightarrow} \varphi_i(U_i) \subset \mathbb{R}^2 \times \mathbb{R}_{\geq 0}$ is a diffeomorphism $(\mathbb{R}_{\geq 0} := \{ t \in \mathbb{R} | t \geq 0\})$ such that 
 if $U_i \cap L_a$ is non-empty, $\varphi_i(U_i \cap L_a) = {\rm Int}(D) \times \{c_a\}$, where  ${\rm Int}(D)$ is the interior of a $2$-disc $D \subset \mathbb{R}^2$ and $c_a \in \mathbb{R}_{\geq 0}$.\\
Here each $L_a$ is called a {\em leaf} of the foliation ${\cal F}$. We call the pair $(M,{\cal F})$ simply a {\em foliated $3$-manifold}.

We call a foliation  ${\cal F}$ on $M$ a {\em complex foliation by Riemann surfaces} and the pair $(M,{\cal F})$ a {\em foliated $3$-manifold by Riemann surfaces} if we require further the following conditions:\\
(3)   each leaf $L_a$ is a Riemann surface, namely, $L_a$ has a complex structure.\\
(4)   the above condition (2)  with replacing  $\mathbb{R}^2$ (leaf coordinate) by $\mathbb{C}$.  \\
(5) $(\varphi_j \circ \varphi_i^{-1})(z,t) = (f_{ij}(z,t), g_{ij}(t))$ for $(z, x) \in \varphi_i(U_i \cap U_j) \subset \mathbb{C} \times \mathbb{R}_{\geq 0}$, where $f_{ij}$ is holomorphic in $z$ and $g_{ij}$ is independent of $z$.\\
\\
{\bf Remark 1.2.} (1) By Frobenius' theorem, the following notions are equivalent if $M$ is closed (cf. [CC; 1.3], [Ta; $\S$ 28]):\\
$\cdot$ $2$-dimensional smooth foliation on $M$,  \\
$\cdot$ $2$-plane field distribution $E$ (rank $2$ subbundle of the tangent bundle $TM$) which is involutive, \\
$\cdot$ $2$-plane field distribution $E$ which is completely integrable. \\
(2)  If the tangent bundle to a $2$-dimensional smooth foliation ${\cal F}$ is orientable, ${\cal F}$ admits a structure of a complex foliation (cf. [HM; Remark 1.2], [Hu; Lemma A.3.1]).  In fact, each orientable $2$-dimensional leaf is equipped with a Riemannian metric and a smooth almost complex structure. By the parametrized Newlander-Nirenberg's theorem ([NN]), it follows that there is a system of complex foliated local coordinate as in the conditions (4) and (5) in Definition 1.1. \\
\\
{\bf Definition 1.3.} For foliated $3$-manifolds $(M,{\cal F})$ and $(M',{\cal F}')$, a {\em morphism} $(M,{\cal F}) \rightarrow (M',{\cal F}')$ is defined to be a smooth map $f : M \rightarrow M'$ such that $f$ maps any leaf $L$ of ${\cal F}$ into a leaf $L'$ of ${\cal F}'$. When $(M,{\cal F})$ and $(M',{\cal F}')$ are foliated $3$-manifolds by Riemann surfaces, we require further $f|_L : L \rightarrow L'$ to be holomorphic for any leaf $L$ of ${\cal F}$. In terms of local coordinates, this is equivalent to saying that for any $p \in M$, there are foliated local coordinates $(U, \varphi)$ and $(U',\varphi')$ around $p$ and $f(p)$, respectively, such that $
(\varphi' \circ f|_U \circ \varphi^{-1})(z,t) = (f^t(z),t)$, where $f^t(z)$ is holomorphic for fixed $t$. 

An {\em isomorphism} $(M,{\cal F}) \stackrel{\sim}{\rightarrow} (M',{\cal F}')$ of foliated $3$-manifolds (resp. foliated $3$-manifolds by Riemann surfaces) is a diffeomorphism $f : M \rightarrow M'$ such that $f|_L : L \rightarrow L'$ is diffeomorphic (resp. biholomorphic) for any leaf $L$ of ${\cal F}$. We identify $(M,{\cal F})$ with $(M',{\cal F}')$ if there is an isomorphism between them.\\
\\
{\bf Definition 1.4.} A {\em smooth dynamical system} (or {\em smooth flow}) on $M$ is defined by a smooth action $\phi$ of $\mathbb{R}$ on $M$, namely, a smooth map $\phi : \mathbb{R} \times M \rightarrow M$ such that $\phi^t := \phi(t, \, \cdot\,)$ is a diffeomorphism of $M$ for each $t \in \mathbb{R}$ which satisfies $\phi^0 = {\rm id}_M, \phi^{t+t'} = \phi^t \circ \phi^{t'}$ for any $t, t' \in \mathbb{R}$.
\\
\\Now we introduce the main object in this paper (cf. [Ko2]).\\
\\
{\bf Definition 1.5.} We define a $3$-dimensional {\em foliated dynamical system} by a triple $\mathfrak{S} = (M, {\cal F}, \phi)$,  
where \\
(1) $M$ is a connected,  closed smooth $3$-manifold,\\
(2) ${\cal F}$ is a complex foliation by Riemann surfaces on $M$,\\
(3) $\phi$ is a smooth dynamical system on $M$,\\
 and these data must satisfy the following conditions (i), (ii): \\
(i) there are finite number of compact leaves $L_1^{\infty}, \dots , L_r^{\infty}$, which may be empty,  such that for any $i$ and $t$ we have $\phi^t(L_i^{\infty}) = L_i^{\infty}$  and that any orbit of the flow $\phi$ is transverse to leaves in $M \setminus \cup_{i=1}^r L_i^{\infty}$. \\
(ii) for each $t \in \mathbb{R}$ the diffeomorphism $\phi^t$ of $M$ maps any leaf to a leaf. \\
In this paper, a $3$-dimensional foliated dynamical system is called an {\em FDS}$^3$ for short.\\
\\
Throughout this paper, we shall use the following notations. For an FDS$^3$ $\mathfrak{S} = (M, {\cal F}, \phi)$, we set
$$\begin{array}{l}
{\cal F}^c := \mbox{ the set of all compact leaves.}\\
{\cal P}_{\mathfrak{S}}^{\infty} := \mbox{the set of non-transverse compact leaves} \\
 \;\; \;\;\;\; \; = \{ L_1^{\infty}, \dots, L_r^{\infty} \} \mbox{ (this set may be empty).}\\
L^{\infty} := \cup_{i=1}^r L_i^{\infty}.\\
M_0 := M \setminus  L^{\infty}. \\
{\cal P}_{\mathfrak{S}} := \mbox{the set of closed orbits of the flow} \; \phi \; \mbox{which are transverse}\\
\;\; \;\;\;\; \;  \;\;\;  \mbox{to leaves in}\;  M_0.  \\
\overline{{\cal P}_{\mathfrak{S}}} := {\cal P}_{\mathfrak{S}} \sqcup {\cal P}_{\mathfrak{S}}^{\infty}.\\
\end{array}
\leqno{(1.6)}
$$
\\
{\bf Remark 1.7.} As Deninger and Kopei suggested (cf. [D1] $\sim$ [D6], [Ko1], [Ko2]), an FDS$^3$ $\mathfrak{S} = (M, {\cal F}, \phi)$  may be regarded as a geometric analogue of a compact arithmetic curve, namely, a smooth proper algebraic curve over a finite field or $\overline{{\rm Spec}({\cal O}_k)} = {\rm Spec}({\cal O}_k) \cup {\cal P}_k^{\infty}$ for the ring ${\cal O}_k$ of integers of a number field $k$, where ${\cal P}_{k}^{\infty}$ is the set of infinite primes of $k$. The set ${\cal P}_{\mathfrak{S}}$ corresponds to the set of closed points (finite primes) of $C$ or ${\rm Spec}({\cal O}_k)$. When  ${\cal P}_{\mathfrak{S}}^{\infty}$ is non-empty, it corresponds to ${\cal P}_k^{\infty}$.  So the analogy is closer if ${\cal P}_{\mathfrak{S}}$ is a countably infinite set. We give such examples in the section 3. \\ 
\\
{\bf Definition 1.8.} Let $\mathfrak{S} = (M, {\cal F}, \phi)$ and $\mathfrak{S}' = (M', {\cal F}', \phi')$ be FDS$^3$'s. An {\em FDS$^3$-morphism} $\mathfrak{S} \rightarrow \mathfrak{S}'$ is defined by a morphism of foliated $3$-manifolds by Riemann surfaces (cf. Definition 1.3) such that  $f : M \rightarrow M'$ commutes with the flows, namely, $f \circ \phi^t = \phi'^t \circ f$ for any $t \in \mathbb{R}$.

An {\em FDS$^3$-isomorphism} $\mathfrak{S} \stackrel{\sim}{\rightarrow} \mathfrak{S}'$ is an FDS$^3$-morphism which is an isomorphism of foliated $3$-manifolds by Riemann surfaces.  We identify $\mathfrak{S}$ with $\mathfrak{S}'$ if there is an FDS$^3$-isomorphism between them.
\\
\\
 {\bf Definition 1.9.} Let $\mathfrak{S} = (M, {\cal F}, \phi)$ be an FDS$^3$. An {\em FDS$^3$-meromorphic function} on $\mathfrak{S}$ is defined to be a smooth map $f : M_0   \rightarrow \mathbb{P}^1(\mathbb{C}) = \mathbb{C} \cup \{ \infty\}$ satisfying the following conditions:\\
 (1) $f$ restricted to any leaf is a meromorphic function,\\
 (2) the zeros and poles of $f$ lie along finitely many closed orbits.\\
\\
For an FDS$^3$ $\mathfrak{S} =(M, {\cal F}, \phi )$, let $T{\cal F}$ denote the subbundle of the tangent bundle $T M_0$  whose total space is the union of the tangent spaces of leaves, and let  $\dot{\phi}^t = \frac{d}{dt} \phi^t$ be the vector field on $M_0$ which generates the flow $\phi$.\\
\\
{\bf Lemma 1.10.} {\em For an FDS$^3$ $\mathfrak{S} =(M, {\cal F}, \phi )$, there is the unique closed smooth $1$-form $\omega$ on $M_0$ satisfying}
$$ \omega|_{T{\cal F}} = 0, \;\; \;\;  \omega(\dot{\phi}^t) = 1. \leqno{({\rm C})}$$
{\em More precisely, let $(M,{\cal F},\phi)$ be a triple satisfying (1), (2), (3) and (i) in Definition 1.5. Then there is the unique smooth $1$-form $\omega$ on $M_0$ satisfying} (C), {\em and the condition} (ii) {\em is equivalent to that $\omega$ is closed.}\\
\\
{\em Proof.} Let $(M,{\cal F},\phi)$ be a triple satisfying (1), (2), (3) and (i) in Definition 1.5. Then we can take a smooth $1$-form $\omega$ on $M_0$ such that for each foliated coordinate
 $(z,t) \in U$, $\omega|_{U} = h(z,t)dt$. Writing $\dot{\phi}^t = a(z,t)\partial_z + b(z,t)\partial_{\overline{z}} + c(z,t)\partial_t$, we have $\omega(\dot{\phi}^t) = c(z,t)h(z,t)$ and so $h(z,t)$ is uniquely determined by $\omega(\dot{\phi}^t) = 1$. This yields the first assertion.

Now, assume the condition (ii). For each point $p \in M_0$, there is a foliated local coordinate $(z,t) \in U$ such that $\phi^s(z,t) = (z, t+s)$. So there is a system of foliated local coordinates $(U_i; \varphi_i)_{i \in I}$ of $M_0$ such that if $U_i \cap U_j$ is not empty,  $(\varphi_j \circ \varphi_i^{-1})(z_i,t_i) = (f_{ij}(z,t), t+c_{ij})$ for some constant $c_{ij} \in \mathbb{R}$. Hence $\omega|_{U_i} = dt_i$ and so $\omega$ is closed.

Conversely, assume that $d\omega = 0$. We may take each $U_i$ is simply connected so that $\omega|_{U_i} = dt_i$ by Poincar\'{e}'s Lemma. Then we can follow the above argument in the reverse way to obtain the condition (ii). $\;\; \Box$\\
\\
{\bf Definition 1.11.} Let $\mathfrak{S} = (M, {\cal F}, \phi)$ be an FDS$^3$. We call the smooth closed $1$-form in Lemma 1.10 the {\em canonical $1$-form} of $\mathfrak{S}$ and denote it by $\omega_{\mathfrak{S}}$. The de Rham cohomology class of  $\omega_{\mathfrak{S}}$ defines the {\em period homomorphism}  
$$ [\omega_{\mathfrak{S}}] : H_1(M_0;\mathbb{Z}) \longrightarrow \mathbb{R}; \;\; [\ell] \mapsto \int_{\ell} \omega_{\mathfrak{S}},$$
and the {\em period group} of $\mathfrak{S}$ is defined by the image of  $[\omega_{\mathfrak{S}}]$, which we denote by $\Lambda_{\mathfrak{S}}$.\\
\\
The following lemma will be used in the examples of the section 3.\\
\\
{\bf Lemma 1.12.} {\em Let $\mathfrak{S}_1 = (M, {\cal F}, \phi_1)$ and $\mathfrak{S}_2 = (M, {\cal F}, \phi_2)$ be FDS$^3$'s. Then we have $\omega_{\mathfrak{S}_1} = \omega_{\mathfrak{S}_2}$ if and only if $(\dot{\phi}_1^t  - \dot{\phi}_2^t)_p \in T_p{\cal F}$ for any $p \in M_0$.}\\
\\
{\em Proof.} Since $\omega_{\mathfrak{S}_i}|_{T{\cal F}} = 0, \omega_{\mathfrak{S}_i}(\dot{\phi}_i^t) = 1$ for $i=1,2$, we have
$$\begin{array}{ll}    (\dot{\phi}_1^t  - \dot{\phi}_2^t)_p \in T_p{\cal F} \;\mbox{for any}\;  p \in M_0 & \Leftrightarrow \omega_{\mathfrak{S}_1}(\dot{\phi}_1^t  - \dot{\phi}_2^t) = 0\\
                                                                                                                                                       &    \Leftrightarrow \omega_{\mathfrak{S}_1}(\dot{\phi}_2^t) =  \omega_{\mathfrak{S}_1}(\dot{\phi}_1^t) = 1\\
                                                                                                                                                       &  \Leftrightarrow \omega_{\mathfrak{S}_1}|_{T{\cal F}} = 0, \omega_{\mathfrak{S}_1}(\dot{\phi}_2^t) = 1\\
                                                                                                                                                       & \Leftrightarrow \omega_{\mathfrak{S}_2} = \omega_{\mathfrak{S}_1}.  
\;\; \;\;\; \Box  \end{array}  
$$\\

 \begin{center}
{\bf 2. A decomposition theorem and a classification of  FDS$^3$'s}
\end{center}

In this section, we give a decomposition theorem for an FDS$^3$, which yields a classification of  FDS$^3$'s. We start with some preparations about the holonomy of a leaf and Tischler's theorem. \\
\\
{\bf 2.1. Foliations without holonomy and Tischler's theorem.} First, we recall the notion of the holonomy group of a leaf (cf. [CC; 2.2, 2.3], [Ta; $\S 22$]). Let $(M, {\cal F})$ be a foliated $3$-manifold and let $L$ be a leaf of ${\cal F}$.  Choose $p \in L$ and let $[c] \in \pi_1(L,p)$ be represented by a loop $c : [0,1] \rightarrow L$ with $c(0)=c(1)=p$. We may choose a subdivision $0 = t_0 < t_1 < \cdots < t_{m+1} = 1$ and a chain of foliated charts  $ (U_0; \varphi_0), \dots, (U_m; \varphi_m)$ such that $c([t_i,t_{i+1}]) \subset U_i$ for $0 \leq i \leq m$. Then  we have 
$$(\varphi_{i+1} \circ \varphi_i^{-1})(z,t) = (f_{i,i+1}(z,t), g_{i,i+1}(t))$$
 for some smooth functions $f_{i,i+1}, g_{i, i+1}$. We set 
$$h_c := g_{m-1,m}\circ \cdots \circ g_{0,1},$$ 
which is a local homeomorphism of $\mathbb{R}$ fixing $0$.  Let $G$ be the group of germs of local homeomorphisms of $\mathbb{R}$ fixing $0$ and let $\widehat{h_c} \in G$ be the germ of $h_c$.  The the correspondence $[c] \mapsto \widehat{h_c}$ gives a well-defined homomorphism
$\Psi : \pi_1(L,p) \rightarrow G$. The {\em holonomy group} ${\cal H}(L)$ of $L$ is defined by the image of $\Psi$, which is independent, up to conjugation, of all choices. A foliated $3$-manifold $(M,{\cal F})$ is said to be {\em without holonomy} if the holonomy group ${\cal H}(L)$ is trivial for any leaf  $L$ of ${\cal F}$.\\
\\
The next theorem, due to Candel and Conlon, is an important ingredient for our decomposition theorem given in the following subsection 2.2. For the proof, we refer to [CC; 9.1]. For a foliated $3$-manifold $(M,{\cal F})$, a subset $X$ of $M$ is said to be {\em ${\cal F}$-saturated} if $X$ is a union of leaves of ${\cal F}$.\\
\\
{\bf Theorem 2.1.1} ([CC; Theorem 9.1.4]). {\em Let $(M, {\cal F})$ be a foliated $3$-manifold with $M$ being closed. Let $X$ be a connected ${\cal F}$-saturated open subset of $M$. Assume that the foliated $3$-manifold $(X,{\cal F}|_X)$ is without holonomy. Then one of the followings holds.}\\
(1) {\em $X$ is a surface bundle over $S^1$ or an open interval, and ${\cal F}|_X$ is the bundle foliation.}\\
(2) {\em Each leaf of ${\cal F}|_X$ is dense in $X$.}\\
\\
Next, we recall Tischler's theorem ([Ti]) and its generalization ([CtC]). \\
\\
{\bf Theorem 2.1.2} ([Ti; Theorem 1]), [CtC; Theorem 2.1]). {\em Let $(M,{\cal F})$ be a foliated $3$-manifold with $M$ being closed. Let $X$ be an ${\cal F}$-saturated open set of $M$. Assume that  there is a non-vanishing, closed $1$-form $\omega$ on $X$ such that  ${\rm Ker}(\omega)$ defines ${\cal F}|_X$. Then $X$ is a surface bundle over $S^1$. }\\
\\
{\bf Remark 2.1.3.}  A non-vanishing, closed $1$-form assumed in Theorem  2.1.2 is approximated (in the $C^{\infty}$-topology) by the $1$-form $\varpi^*(d\theta)$, where $\varpi : U \rightarrow S^1$ is the fibration and $d\theta$ is an angular $1$-form on $S^1$ (cf. ibid). \\
\\
{\bf 2.2. A decomposition theorem.}  Let $\mathfrak{S} =(M, {\cal F}, \phi )$ be an FDS$^3$. Let $X_{1}, \dots , X_d$ be connected components of $M_0 := M \setminus L^{\infty}$.  Each $(X_a, {\cal F}|_{X_a})$ is a foliated $3$-manifold for $1\leq a \leq d$. \\
\\
{\bf Lemma 2.2.1.} {\em The foliated 3-manifold $(X_a, {\cal F}|_{X_a})$ is without holonomy.}\\
\\
{\em Proof.} Let $L$ be any leaf of $ {\cal F}|_{X_a}$. Since $\phi^t$ maps any leaf to a leaf, there is a system foliated charts $\{ (U_i; \varphi_i)\}$ which covers $L$ and satisfies
$$  \varphi_i(\phi^t(p)) = (z_i(p), t)$$
for any $p \in L$. Therefore, if $U_i \cap U_j$ is non-empty, we have
$$ (\varphi_i \circ \varphi_j^{-1})(z_j,t) = (z_i(\varphi_j^{-1}(z_j,t)), t)$$
and so the holonomy group ${\cal H}(L)$ is trivial. Hence $(X_a, {\cal F}|_{X_a})$ is without holonomy. $\;\; \Box$\\
\\
{\bf Theorem 2.2.2} ({\em A decomposition theorem}). {\em Let $\mathfrak{S} = (M, {\cal F}, \phi)$ be an FDS$^3$ and let $X_{a}$ be a connected component of $M_0 := M \setminus L^{\infty}$. The foliated 3-manifold $(X_a, {\cal F}|_{X_a})$ is one of the followings.}\\
(1) {\em $X_a$ is a surface bundle over $S^1$ or an open interval, and ${\cal F}|_{X_a}$ is the bundle foliation.}\\
(2)  {\em $X_a$ is a surface-bundle over $S^1$ and any leaf in ${\cal F}|_{X_a}$ is dense in $X_a$.}\\
\\
{\em Proof.}  By Definition 1.5 of an FDS$^3$, $X_a$ is a connected ${\cal F}$-saturated open subset of $M$. By Theorem 2.1.1, (1) $X_a$ is a surface bundle over $S^1$ or an open interval, and ${\cal F}|_X$ is the bundle foliation, or (2) each leaf of ${\cal F}|_X$ is dense in $X$. For the case (2), $X$ is a surface bundle over $S^1$ by Lemma 1.10 and Theorem 2.1.2. $\;\; \Box$\\
\\
{\bf Remark 2.2.3.} For an FDS$^3$ $\mathfrak{S} =(M, {\cal F}, \phi )$, if $M$ is cut along $L^{\infty}$, then $M_0$ is decomposed into connected components $M_0 = X_1 \sqcup \cdots \sqcup X_d$, where each foliated $3$-manifold $(X_a, {\cal F}|_{X_a})$ has the structure described in Theorem 2.2.2 and the flow $\phi|_{X_a}$ is transverse to leaves of ${\cal F}|_{X_a}$. This is the reason that we call Theorem 2.2.2 a {\em decomposition theorem}. It may remind us of the JSJ decomposition of a $3$-manifold ([JS], [Jo]).\\
\\
Theorem 2.2.2 can be restated as the following classification of FDS$^3$'s. For the notations, see (1.6).\\
\\
{\bf Corollary 2.2.4} ({\em A classification}). {\em An FDS$^3$ $\mathfrak{S} = (M, {\cal F}, \phi)$ is classed as one of the following  types}:\\
I. {\em ${\cal F} = {\cal F}^c$ and ${\cal P}_{\mathfrak{S}}^{\infty}$ is empty}. \\
@{\em Then $M$ is a surface bundle over  $S^1$ and ${\cal F}$ is the bundle foliation.}\\
II. {\em ${\cal F}^c$ is empty}.  \\
@{\em Then $M$ is a surface bundle over $S^1$ and any leaf of ${\cal F}$ is dense in $M$.}\\
III. {\em ${\cal P}_{\mathfrak{S}}^{\infty}$ is a non-empty $($finite$)$ set. Let $X_{a}$ be a connected component of $M_0$.} \\
@{\em Then the foliated 3-manifold $(X_a, {\cal F}|_{X_a})$ is one of the followings.}\\
III-1. {\em $X_a$ is a surface bundle over $S^1$ and ${\cal F}|_{X_a}$ is the bundle foliation.}\\
III-2. {\em $X_a$ is a surface bundle over an open interval and ${\cal F}|_{X_a}$ is the bundle foliation.}\\
III-3. {\em $X_a$ is a surface-bundle over $S^1$ and any leaf in ${\cal F}|_{X_a}$ is dense in $X_a$.}\\
\\
\\
We note that the class of FDS$^3$'s of bundle foliation over $S^1$ is characterized by the period group. Although this may be known  (cf. [CCI; 9.3], [Fa; 2.1]), we give a proof in the following, for the sake of readers.\\
\\ 
{\bf Proposition 2.2.5.}  {\em Let $\mathfrak{S}$ be an FDS$^3$. If $\mathfrak{S}$ is of type I or of type III-1, then  the period group $\Lambda_{\mathfrak{S}} = \mathbb{Z}$. Conversely, if the $\Lambda_{\mathfrak{S}}$ has rank one (namely,  $\Lambda_{\mathfrak{S}} \simeq \mathbb{Z}$) and $M_0$ is connected, then $\mathfrak{S} $ is of type I or of type III-1.}\\
\\
{\em Proof.} Suppose that $\mathfrak{S} =  (M, {\cal F},\phi)$ is an FDS$^3$ of type I or of type III-1.  Then, for any connected component $X_a$ of $M_0$, there is a fibration $\varpi_a : X_a \rightarrow S^1$.  Let $d\theta$ be the angular $1$-form on $S^1$ such that $\int_{S^1} d\theta = 1$. Then the canonical $1$-form $\omega_{\mathfrak{S}}$ is given by $\omega_{\mathfrak{S}}|_{X_a} = \varpi_a^{*}(d\theta)$ for any $a$.  For $[\ell] \in H_1(X_a;\mathbb{Z})$, we have
$$ [\omega_{\mathfrak{S}}]([\ell])  = \int_{(\varpi_a)_*(\ell)} d\theta = \; \mbox{\rm the degree of}\; \varpi_*(\ell) \; \mbox{\rm on}\; S^1$$
for any $a$ and hence $\Lambda_{\mathfrak{S}} = \mathbb{Z}$. Conversely, suppose that $\mathfrak{S} =  (M, {\cal F},\phi)$ is an FDS$^3$ with $M_0$ being connected and $\Lambda_{\mathfrak{S}} = \lambda\mathbb{Z}$ for some $\lambda \in \mathbb{R}^{\times}$ so that $ [\omega_{\mathfrak{S}}](H_1(M_0;\mathbb{Z})) = \lambda \mathbb{Z}.$ Fix a base point $p_0 \in M_0$ and define the map $\varpi : M_0 \rightarrow S^1$ by
 $$ \varpi(p) := \exp \left( \frac{2\pi i}{\lambda} \int_{\gamma} \omega_{\mathfrak{S}}\right),$$
 where $\gamma$ is a path from $p_0$ to $p$. Then we see easily that $\lambda \varpi^{*}(d\theta) = \omega_{\mathfrak{S}}$. By Definition 1.11 of $\omega_{\mathfrak{S}}$, $\varpi$ is a fibration and ${\cal F}$ consists of fibers of $\varphi$. Hence $\mathfrak{S}$ is of type I or of type III-1. $\;\; \Box$\\
\\
Finally we note that an FDS$^3$ of type III-2 has no transverse closed orbits.\\
\\
{\bf Proposition 2.2.6.} {\em  Let $\mathfrak{S} = (M,{\cal F},\phi)$ be an FDS$^3$ of type III-2. Then ${\cal P}_{\mathfrak{S}}$ is empty and $\Lambda_{\mathfrak{S}} = \{ 0 \}$.}\\
\\
{\em Proof.}  Let $X_a$ be a connected component of $M_0$. We may identify $(X_a, {\cal F}|_{X_a})$ with $(L \times (0,1), \{ L \times \{t \} \}_{t \in (0,1)})$, where $L \in {\cal F}|_{X_a}$. Let $\varpi : X_a \rightarrow (0,1)$ be
the projection. For any   closed curve $\ell : S^1 \rightarrow X_a$, $\varpi \circ \ell : S^1 \rightarrow (0,1)$ has the maximum $\mu$ and so $\ell$ is not transverse to the leaf
$L \times \{\mu \}$. Hence ${\cal P}_{\mathfrak{S}}$ is empty. Next, since $X_a$ is homotopy equivalent to the leaf $L \times \{ \frac{1}{2} \}$,
$\ell$ is homotopic to a closed curve $\ell'$ in $L \times \{ \frac{1}{2} \}$. Since $\omega_{\mathfrak{S}}|_{T{\cal F}} = 0$, we have $[\omega_{\mathfrak{S}}]([\ell]) = [\omega_{\mathfrak{S}}]([\ell']) = 0$. Hence $\Lambda_{\mathfrak{S}} = \{ 0 \}$.
$\;\; \Box$\\

\begin{center}
{\bf 3.  Examples of FDS$^3$'s}
\end{center}

In this section, we construct concrete examples of FDS$^3$'s for each type in Corollary 2.2.4. \\
\\
I. We give an example of an FDS$^3$ of type I. Note that any smooth surface bundle over $S^1$ is obtained by the mapping torus of a surface diffeomorphism.\\
\\
{\bf Example 3.I} (Mapping torus and pseudo-Anosov flow).  Let $\Sigma_g$ be a connected, closed smooth surface of genus $g \geq 1$ and let $\varphi$ be a diffeomorphism of $\Sigma_g$. 
Let $M$ be the {\em mapping torus}  $M(\Sigma_g, \varphi)$ defined by 
$$ M(\Sigma_g, \varphi) := (\Sigma_g \times [0, 1])/(z, 1) \sim (\varphi(z), 0).$$
Then the projection
$$\varpi : M \to S^1 = \mathbb{R}/\mathbb{Z}; \;\; \varpi([z, s]) = s \; (\mbox{mod} \; \mathbb{Z}) $$
is a fibration, where each fiber $\varpi^{-1}(\theta)$ over $\theta \in S^1$ is diffeomorphic to $\Sigma_g$. The set of fibers ${\cal F} := \{ \varpi^{-1}(\theta) \}_{\theta \in S^1}$ defines a $2$-dimensional foliation, the bundle foliation.  Since $T{\cal F}$ is orientable, ${\cal F}$ admits a complex foliation structure (cf. Remark 1.2 (2)). Let  $\phi$ be the suspension flow defined by 
$$\phi^t([z, s]) := [z, s+t].$$
Then $\mathfrak{S} := (M, {\cal F},\phi)$ forms an FDS$^3$. Suppose further that the diffeomorphism $\varphi$ is of pseudo-Anosov type (cf. [Cl; 1.11]). Note that $\varphi$ is an Anosov diffeomorphism when $\Sigma_g$ is a torus. Then, since $\varphi$ has countably infinite periodic points,  $P_{\mathfrak{S}}$ is a countably infinite set. By Proposition 2.2.5, the canonical 1-form $\omega_{\mathfrak{S}}$ is $\varpi^*(d\theta)$ and the period group $\Lambda_{\mathfrak{S}}  =  \mathbb{Z}$.
\\
 \\
{\bf Remark 3.1.} An FDS$^3$ of type I may be regarded as an analogue of a smooth proper algebraic curve $C$ over a finite field $\mathbb{F}_q$, where the $2$-dimensional foliation corresponds to  the geometric fiber $C \otimes \overline{\mathbb{F}_q}$ and the monodromy $\varphi$  corresponds to the Frobenius automorphism in ${\rm Gal}(\overline{\mathbb{F}_q}/\mathbb{F}_q)$. \\
\\
 II. We construct two examples of FDS$^3$'s of type II. First, we give a method to construct countably infinitely many closed orbits using the horseshoe map. 
 Let $\mathfrak{S} = (M, {\cal F}, \phi)$ be an FDS$^3$ and $\gamma \in {\cal P}_{\mathfrak{S}}$. We say that the flow $\phi$ is {\em of contraction type} around $\gamma$ if there is 
 a tubular neighborhood $V = D \times S^1$ of $\gamma$, where $D (\subset \mathbb{C})$ is a $2$-disc centered at $0$ and $\gamma = \{ 0\} \times S^1$, such that for any $t \in \mathbb{R}$ and $p \in \gamma$, $\phi^t(V) \subset V$ and 
$$ \phi^t_{D,p} := \phi^t|_{D \times \{ p\}} : D = D \times \{p\} \stackrel{\phi^t}{\longrightarrow}  \phi^t(D) \times \{ \phi^t(p) \} \subset D$$
is a contraction map, namely, $T_z(\phi^t_{D,p})$ has the eigenvalues $\lambda_1, \lambda_2$  satisfying $0 < |\lambda_1|, |\lambda_2| < 1$ for $z \in {\rm Int}(D) \setminus \{ 0\}$. Let $h : {\rm Int}(D) \rightarrow {\rm Int}(D)$ be the {\em horseshoe diffeomorphism} ([Sm]) and let $U = M({\rm Int}(D), h)$ be the mapping torus of $h$ equipped with the suspension flow $\phi_h$. Note that there is a diffeomorphism  $\psi : U \stackrel{\approx}{\rightarrow} {\rm Int}(V)$. 
\\
\\
{\bf Lemma 3.2.} {\em Notations being as above, let us replace $({\rm Int}(V), \phi|_{{\rm Int}(V)})$ by $(U, \phi_h)$ via $\psi$ so that the resulting $3$-manifold $M_{\gamma,h}$ is equipped with a new smooth flow $\phi_{\gamma,h}$ satisfying $\phi_{\gamma,h} = \phi_{h}$ in $U$ and $\phi_{\gamma,h} = \phi$ in $M \setminus {\rm Int}(V)$. We define the  foliation ${\cal F}_{\gamma,h}$ on $M_{\gamma,h}$ by the foliation ${\cal F}$ on $M$ via $\psi$. Then the triple $\mathfrak{S}_{\gamma, h} := (M_{\gamma,h}, {\cal F}, \phi_{\gamma,h})$ forms an FDS$^3$ and  the flow $\phi_{\gamma, h}$ has countably infinitely many closed orbits around $\gamma$. Moreover, the canonical forms of $\mathfrak{S}$ and $\mathfrak{S}_{\gamma, h}$ are same, $\omega_{\mathfrak{S}} = \omega_{\mathfrak{S}_{\gamma, h}}$.}  \\
\\
{\em Proof.} That $(M_{\gamma,h}, {\cal F}_{\gamma, h}, \phi_{\gamma,h})$ forms an FDS$^3$ is easily seen by the construction. The latter property follows from that the horseshoe map has countably infinitely many periodic points ([KH; Cororally 2.5.1]). For any $p \in \gamma$, we may assume, by the change of parameters, that there is $t_0 > 0$ such that $\phi^{t_0}|_{{\rm Int}(D) \times \{ p\}} : {\rm Int}(D) \times \{p\} \rightarrow \phi^{t_0}({\rm Int}(D)) \times \{ p\}$ and $\phi_h^{t_0}|_{{\rm Int}(D) \times \{ p\}} : {\rm Int}(D) \times \{p\} \rightarrow \phi_h^{t_0}({\rm Int}(D)) \times \{ p\}$ are first return diffeomorphisms. Since these are isotopic and an isotopy is realized as a flow in ${\rm Int}(V)$, we see that the canonical $1$-forms are unchanged.  $\;\; \Box$\\
\\
{\bf Example 3.II.1.} The following example of a foliated $3$-manifold was considered in [Cl; 4.2] and  [Hu; A.5] (see also Alvarez-Lopez's example in [D7; page 10]).  Let $T^2$ be the $2$-dimensional torus $\mathbb{C}/\mathbb{Z}^2$ so that the fundamental group $\pi_1(T^2)$ is generated by the homotopy classes of a meridian and a longitude, say $\mathfrak{m}$ and $\mathfrak{l}$, respectively. Let $\rho : \pi_1(T^2) \rightarrow \mathbb{R}$ be a given homomorphism such that $\rho(\mathfrak{m}) \notin \mathbb{Q}$ or $\rho(\mathfrak{l}) \notin \mathbb{Q}$, and we set $\overline{\rho}(g) := \rho(g) \; \mbox{mod} \; \mathbb{Z} \in S^1$ for $g \in \pi_1(T^2)$. Let $M$ be the quotient $3$-manifold of $\mathbb{C} \times S^1$ by the action of $\pi_1(T^2)$
$$ M := (\mathbb{C} \times S^1)/\pi_1(T^2),$$
where $\pi_1(T^2)$ acts on the universal cover $\mathbb{C}$ of $T^2$ as the monodromy and on $S^1$ by $\theta \mapsto \theta + \overline{\rho}(g)$ for $\theta \in S^1$ and $g \in \pi_1(T^2)$. Let $L_{\theta}$ denote the image in $M$ of $\mathbb{C} \times \{\theta \}$ in $M$. Then ${\cal F} := \{ L_{\theta}\}_{\theta \in S^1}$ forms a $2$-dimensional foliation on $M$. Here we see that
$$ L_\theta = \left\{ 
\begin{array}{ll}
S^1 \times \mathbb{R} &  \; \mbox{if}\; \rho(\mathfrak{m}) \notin \mathbb{Q}, \rho(\mathfrak{l}) \in \mathbb{Q} \; \mbox{or} \; \rho(\mathfrak{m}) \in \mathbb{Q},  \rho(\mathfrak{l}) \notin \mathbb{Q},\\
\mathbb{R}^2 & \; \mbox{if}\; \rho(\mathfrak{m}), \rho(\mathfrak{l}) \notin \mathbb{Q} \; \mbox{and}\; \rho(\mathfrak{m})/\rho(\mathfrak{l}) \notin \mathbb{Q},
\end{array} \right.
$$
and that any leaf $L_\theta$ is dense in $M$. Since $T{\cal F}$ is orientable, ${\cal F}$ admits a complex foliation structure (Remark 1.2 (2)). Let $\phi_1$ be the flow defined by
$$ \phi_1^t([z, \theta])   := [z, \theta + \overline{t}]$$
where $\overline{t} := t$ mod $\mathbb{Z}$. Then we obtain an FDS$^3$ $\mathfrak{S}_1 := (M, {\cal F}, \phi_1)$ of type II. The canonical 1-form $\omega_{\mathfrak{S}_1}$ is given by $\omega_{\mathfrak{S}_1}|_{T{\cal F}} = 0$ and $\omega_{\mathfrak{S}_1}|_{TS^1} = d\theta$, and the period group $\Lambda_{\mathfrak{S}_1} = \mathbb{Z} + \rho(\mathfrak{m})\mathbb{Z} + \rho(\mathfrak{l})\mathbb{Z}$. For this example, however, any orbit of the flow $\phi_1$  is closed with period $1$ and hence $P_{\mathfrak{S}}$ is uncountable.

In order to obtain an FDS$^3$ $\mathfrak{S}$ with countably infinite ${\cal P}_{\mathfrak{S}}$, we interpret the above $(M, {\cal F})$ from a different view and define a new dynamical system. In fact, the natural map $\mathbb{C} \times S^1 \rightarrow T^2 \times S^1$ induces the diffeomorphism
$$ M \; \stackrel{\approx}{\longrightarrow} \; T^3 = (\mathbb{R}/\mathbb{Z})^3.$$
Then the leaf $L_\theta$ is given by
$$L_\theta := \{ (\theta_1, \theta_2, \theta - \rho(\mathfrak{m}) \theta_1 - \rho(\mathfrak{l}) \theta_2) \, | \, \theta_1, \theta_2 \in S^1 \}$$
 for $p = (\theta_1, \theta_2, \theta_3) \in T^3$. 
Let $V_{\mathfrak{m}} := ( 1, 0,  -\rho(\mathfrak{m})), V_{\mathfrak{l}} := (0, 1, -\rho(\mathfrak{l}))$ be vector fields on $T^3$, and 
we define the smooth dynamical system $\phi_2$ on $T^3$ by the equation
$$ \begin{array}{ll} \displaystyle{\frac{d}{dt}\phi_2^t(p)} & :=(0,0,1)+\sin(2\pi\theta_1)V_{\mathfrak{m}}+\sin(2\pi\theta_2)V_{\mathfrak{l}} \\
& =(\sin(2\pi\theta_1), \sin(2\pi\theta_2), 1-\rho(\mathfrak{m})\sin(2\pi\theta_1)-\rho(\mathfrak{l})\sin(2\pi\theta_2)).
\end{array} $$
Then we obtain an FDS$^3$ $\mathfrak{S}_2 := (T^3, {\cal F}, \phi_2)$ of type II. Since $V_{\mathfrak{m}}$ and $V_{\mathfrak{l}}$ are tangent to ${\cal F}$ and so $(\dot{\phi}_1^t - \dot{\phi}_2^t)_p \in T_p{\cal F}$ for $p \in M = T^3$, we have $\omega_{\mathfrak{S}_2} = \omega_{\mathfrak{S}_1}, \Lambda_{\mathfrak{S}_2} = \Lambda_{\mathfrak{S}_1}$ by Lemma 1.12. We see that ${\cal P}_{\mathfrak{S}_2}$ consists of the following four closed orbits
$$ \left\{ \begin{array}{ll}
\gamma_1 = \{ (0,0, \theta) \, |\,  \theta \in S^1 \},\\
\gamma_2 = \{ (\frac{1}{2}, 0, \theta) \, | \, \theta \in S^1 \}, \; \gamma_2' = \{ (0, \frac{1}{2}, \theta) \, | \, \theta \in S^1 \},\\
\gamma_3 = \{  (\frac{1}{2},  \frac{1}{2}, \theta) \, | \, \theta \in  S^1 \},\\
\end{array} \right.
$$
and that $\phi_2$ is of contracting type around  $\gamma_3$. We replace a tubular neighborhood of $\gamma_3$ by the suspension of the horseshoe map. By Lemma 3.2,  we obtain an FDS$^3$ $\mathfrak{S} = (M, {\cal F}, \phi)$ of type II such that  ${\cal P}_{\mathfrak{S}}$ is countably infinite, and the canonical 1-form $\omega_{\mathfrak{S}}$ and the period group $\Lambda_{\mathfrak{S}}$ are same as $\omega_{\mathfrak{S}_i}$ and $\Lambda_{\mathfrak{S}_i}$ ($i=1,2$), respectively.\\
\\
{\bf Example 3.II.2.} Let $T^2 := (\mathbb{R}/\mathbb{Z})^2$ be the $2$-dimensional torus and let $\varphi_A : T^2 \rightarrow T^2$ be the linear Anosov diffeomorphism of $T^2$ defined by the matrix 
$$ A = \left(   \begin{array}{cc} 3 & 1 \\ 2 & 1 \end{array} \right).$$
 Then there are  two fixed points $ p_1 = (0, 0)$ and $p_2 = (\frac{1}{2}, 0)$ of $\varphi_A$.  
For each $p_i$, we remove $p_i$ from $T^2$  and glue $(T_{p_i}(T^2) \setminus \{ 0\})/\mathbb{R}_+ = S^1$ on it. Then we obtain a twice-punctured torus $T^2_*$ with $\partial T^2_* = S^1_1 \sqcup S^1_2$, $S^1_i = S^1$, and the diffeomorphism $ \varphi_A^* : T^2_* \rightarrow T^2_*.$

Let $M^* = M(T^2_*, \varphi_A^*)$ be the mapping torus of $\varphi_A^*$, which  has two boundary components, say $T_1$ and $T_2$, which correspond to $p_1$ and $p_2$, respectively. Let ${\cal F}^* = \{ L_\theta \}_{\theta \in S^1}$ be the bundle foliation of the fibration $\varpi_* : M^* \rightarrow S^1$ and let $\phi_*$ be the suspension flow of $\varphi_A^*$. For each $i$, ${\cal F}^*|_{T_i}$ defines a foliation on $T_i$ whose leaves are $S^1$. 
The $1$-form $\omega_* := \varpi_*^{-1}(d\theta)$ satisfies $\omega_*|_{T{\cal F}^*} = 0, \omega_*(\dot{\phi}_*^t) = 1$.

Finally, let $\lambda \in \mathbb{R} \setminus \mathbb{Q}$ and consider the diffeomorphism $\psi_{\lambda} : T_1 \rightarrow T_2$ defined by sending $L_{\theta}$ to $L_{\theta + \lambda}$. We define $M$ to be the $3$-manifold obtained from $M^*$ by gluing $T_1$ and $T_2$ via $\psi_{\lambda}$. In fact, $M$ is a $\Sigma_2$-bundle over $S^1$. The foliation ${\cal F}$ and the flow $\phi$ on $M$ are induced by 
 ${\cal F}^*$ and $\phi_*$, and ${\cal F}$ has a complex structure. Then any leaf in ${\cal F}$ is dense in $M$ by the way of the gluing. Thus we obtain the FDS$^3$ $\mathfrak{S} = (M, {\cal F}, \phi)$ of type II. Further, since $\phi$ comes from the suspension of the Anosov diffeomorphism, ${\cal P}_{\mathfrak{S}}$ is a countably infinite set. The canonical $1$-form $\omega_{\mathfrak{S}}$ is induced by $\omega_*$ and the period group $\Lambda_{\mathfrak{S}} = \mathbb{Z} + \lambda\mathbb{Z}$.
\\
\\
III. We construct  examples of FDS$^3$'s of type  III. This type of FDS$^3$ may be regarded as an analogue of a number ring $\overline{ {\rm Spec}({\cal O}_k)} = {\rm Spec}({\cal O}_k) \cup {\cal P}_k^{\infty}$ in the respect that ${\cal P}_{\mathfrak{S}}^{\infty}$ corresponds to the set ${\cal P}_k^{\infty}$ of infinite primes of a number field $k$ (cf. Remark 1.7).\\
\\
III-1. We give two examples of FDS$^3$'s of type III-1.\\
\\
{\bf Example 3.III-1.1} (Reeb foliation on $S^3$ and the horseshoe flow).  Let $M$ be the $3$-sphere $S^3$. Let $M := S^3 = V_1 \cup V_2$ be the Heegaard splitting of genus one, where $V_i$ is a solid torus $D^2 \times S^1$ ([He; Chapter 2]). Consider the {\em Reeb foliation} ${\cal F}_i$ on each $V_i$ and so the $2$-dimensional foliation ${\cal F}$ on $S^3$ by getting ${\cal F}_1$ and ${\cal F}_2$ together, where any leaf is diffeomorphic to $\mathbb{R}^2$ besides the only one compact leaf $L^{\infty} = \partial V_1 = \partial V_2$ ([CC; 1.1], [Ta; $\S 1$]). We define the dynamical system $\phi$ as follows.

First, we consider a flow $\phi_1$ on $M$ such that any orbit of $\phi_1$ is transverse to leaves in ${\cal F} \setminus L^{\infty}$ and $\phi_1^t(L^{\infty}) = L^{\infty}$ for $t \in \mathbb{R}$. Then  there is the only one closed orbit $\gamma_i = \{ 0 \} \times S^1$ in each $V_i$ ($0$ being the center of $D^2$). In fact, the flow $\phi_1$ on ${\rm Int}(V_i)$ is the suspension flow of a contraction map $ \varphi_i : {\rm Int}(D^2) \ni z \mapsto az \in {\rm Int}(D^2)$ $(0<a<1)$ and so $\phi_1$ is of contraction type around $\gamma_i$.

Next, we replace the contraction map $\varphi_i$ by the horseshoe map $h$ and  the flow $\phi_1$ around $\gamma_i$ by the suspension of $h$.  Let $\phi$ be the resulting flow on all of $M$. By Lemma 3.2  we have an
 FDS$^3$ $\mathfrak{S} = (M, {\cal F},\phi)$ of type III-1 such that  ${\cal P}_{\mathfrak{S}}$ is a countably infinite set and ${\cal P}_{\mathfrak{S}}^{\infty}$ consists of the only one non-transverse compact leaf $L^{\infty}$. We have the decomposition of $M_0 = S^3 \setminus L^{\infty}$ into connected components 
$$M_0 =  {\rm Int}(V_1) \sqcup {\rm Int}(V_2),$$ 
and each $({\rm Int}(V_i), {\cal F}|_{{\rm Int}(V_i)})$ ($i=1,2$) is the bundle foliation of the ${\rm Int}(D^2)$-bundle over $S^1$. By Proposition 2.2.5, the canonical 1-form $\omega_{\mathfrak{S}}$ restricted on each ${\rm Int}(V_i)$ is the pull-back of  the angular $1$-form of $S^1$ under the fibration, and so the period group $\Lambda_{\mathfrak{S}} = \mathbb{Z}$.

In view of the analogy in arithmetic topology ([Mo; Chapter 3]), the $3$-sphere $S^3$ may be regarded as an analogue of $\overline{ {\rm Spec}(\mathbb{Z})} = {\rm Spec}(\mathbb{Z}) \cup {\cal P}_{\mathbb{Q}}^{\infty}$, where ${\cal P}_{\mathbb{Q}}^{\infty}$ consists of the only one infinite prime of $\mathbb{Q}$. \\
\\
{\bf Example 3.III-1.2} (Open book decomposition). Let $M$ be a closed $3$-manifold. It is known that $M$ contains a fibered link $L = K_1 \cup \cdots \cup K_r$, namely, there is a fibration $\varpi : M \setminus  {\rm Int}(V(L)) \rightarrow S^1$, where ${\rm Int}(V(L)) = \sqcup_{i=1}^r {\rm Int}(V(K_i))$ is the interior of a tubular neighborhood of $L$ and any fiber of $\varpi$ is a surface with $r$ boundary components. We have the foliation on $M \setminus {\rm Int}(V(L))$ by tubularizing the fibers of $\varpi$ around $\partial V(L) = \sqcup_{i=1}^r  \partial {\rm Int}(V(K_i))$ ([CC; Example 3.3.11]). The structure this induces on $M$ is called an {\em open book decomposition} ([Cl; Example 4.11]). We fill in $V(L)$ with the Reeb component to obtain the foliation ${\cal F}$ on all of $M$. We define the flow on $M \setminus {\rm Int}(V(L))$ by the suspension of the monodromy $\varphi$ of the fibration $\varpi$. We suppose that $\varphi$ is of pseudo-Anosov type (for example, this is the case if $L$ is a hyperbolic link). The flow on ${\rm Int}(V(L))$ is defined to be the one transverse to any leaf of the Reeb foliation. Thus we have an FDS$^3$ $\mathfrak{S} = (M, {\cal F},\phi)$ of type III such that  ${\cal P}_{\mathfrak{S}}^{\infty} = \{ \partial V(K_1), \dots , \partial V(K_r)  \}$ and ${\cal P}_{\mathfrak{S}}$ is a countably infinite set. We have the decomposition of $M_0 = M \setminus \partial V(L)$ into connected components 
$$M_0 = \sqcup_{i=1}^r {\rm Int}(V(K_i)) \sqcup (M \setminus V(L)).$$
 Here ${\rm Int}(V(K_i)) = \mathbb{R}^2 \times S^1$ and ${\cal F}|_{{\rm Int}(V(K_i))}$ is the bundle foliation over $S^1$, and $M \setminus V(L)$ is a surface bundle over $S^1$ and ${\cal F}|_{M \setminus V(L)}$ is also the bundle foliation over $S^1$. By Proposition 2.2.5, the canonical 1-form $\omega_{\mathfrak{S}}$ restricted on ${\rm Int}(V(K_i))$ or $M \setminus V(L)$ is the pull-back of  the angular $1$-form of $S^1$ under the fibration, and so the period group $\Lambda_{\mathfrak{S}} = \mathbb{Z}$.\\
\\
{\bf Remark 3.3.} Example 3.III-1.2 shows that any closed smooth $3$-manifold $M$ admits a structure of an FDS$^3$ with non-empty $P_{\mathfrak{S}}^{\infty}$. So the notion of an FDS$^3$ is generic in this sense. Moreover,   the fibration $\varpi : M \setminus   {\rm Int}(V(L)) \rightarrow S^1$ above induces the surjective homomorphism $\pi_1(M \setminus L) \rightarrow \mathbb{Z}$ and hence $M$ has a $\mathbb{Z}$-covering ramified over $L$. This may be analogous to that any number field $k$ has a $\mathbb{Z}_l$-extension ramified over primes $\mathfrak{l}_1, \dots, \mathfrak{l}_n$ lying above $(l)$, where $l$ is a prime number. \\
\\
III-2. We give an example of an FDS$^3$ of type III-2.\\
\\
{\bf Example 3.III-2.} Let $T^3 := (\mathbb{R}/\mathbb{Z})^3 =\{ (\theta_1, \theta_2, \theta_3) \, | \, \theta_i  \in \mathbb{R}/\mathbb{Z} \}$ be the $3$-dimensional torus and let $\mathcal{F}=\{ T^2 \times \{ \theta_3 \} \}_{\theta_3 \in \mathbb{R}/\mathbb{Z}}$ be the linear foliation which admits a complex structure.
We define the smooth dynamical system $\phi$  on $T^3$ by 
$$ \frac{d}{dt} \phi^t(p) := ( \cos(2\pi \theta_3), 0, \sin(2\pi \theta_3) ) $$
for $p = (\theta_1, \theta_2, \theta_3)$. Then we obtain an FDS$^3$ $(T^3, \mathcal{F}, \phi)$ which is equipped with ${\cal P}^{\infty}_{\mathfrak{S}} = \{ L_1^{\infty}, L_2^{\infty}\}$, where 
$$L_1^{\infty} :=\{(\theta_1, \theta_2, 0) \, | \, \theta_i, \in  \mathbb{R}/ \mathbb{Z} \}, L_2^{\infty} :=\{ (\theta_1, \theta_2 ,\frac{1}{2}) \, | \, \theta_i  \in \mathbb{R}/ \mathbb{Z}  \}.$$
We have the decomposition of $M_0 := T^3 \setminus (L_1^{\infty} \cup L_2^{\infty}) $  into connected components 
$$ \begin{array}{l} M_0 = X_1 \sqcup X_2, \\
X_1 :=\{ (\theta_1, \theta_2, \theta_3) \, | \, \theta_1, \theta_2 \in \mathbb{R}/\mathbb{Z}, 0<\theta_3<\frac{1}{2}\}, \\
X_2 :=\{(\theta_1, \theta_2, \theta_3) \, | \, \theta_1, \theta_2 \in \mathbb{R}/\mathbb{Z}, \frac{1}{2}<\theta_3<1\}.
\end{array}$$
Here each $X_a$ $(a=1,2)$ is  a $T^2$-bundle over an open interval and ${\cal F}|_{X_a}$ is the bundle foliation.  The canonical $1$-form $\omega_{\mathfrak{S}}$ is given by $\omega_{\mathfrak{S}}|_{X_a} = {\rm cosec}(2\pi \theta_3) d\theta_3$ and the period group $\Lambda_{\mathfrak{S}} =\{0\}$ by Proposition 2.2.6D\\
\\
III-3. Finally, we give an example of an FDS$^3$ of type III-3.\\
\\
{\bf Example 3.III-3.} Let $T^3 := (\mathbb{R}/\mathbb{Z})^3 =\{ (\theta_1, \theta_2, \theta_3) \, | \, \theta_i  \in \mathbb{R}/\mathbb{Z} \}.$ Fix an irrational number $\rho \in \mathbb{R} \setminus \mathbb{Q}$. 
Let $\omega_0$ be the smooth $1$-form on $T^3$ defined by 
$$\omega_0 :=\sin(2\pi \theta_3)(d\theta_1+\rho d\theta_2)+d\theta_3.$$
Since we see $\omega_0 \wedge d\omega_0 =0$, ${\rm Ker}(\omega_0)$ defines the foliation $\mathcal{F}$ on $T^3$ by Frobenius' theorem (cf. Remark 1.2 (1)). 
Let $\phi$ be the smooth dynamical system defined by
$$ \frac{d}{dt} \phi_1^t(p) := (1, 0, 0) $$
for $p = (\theta_1, \theta_2, \theta_3)$. Then we obtain an FDS$^3$ $\mathfrak{S}_1 := (T^3, \mathcal{F}, \phi_1)$ which is equipped with ${\cal P}^{\infty}_{\mathfrak{S}} = \{ L_1^{\infty}, L_2^{\infty}\}$, where 
$$L_1^{\infty} :=\{(\theta_1, \theta_2, 0) \, | \, \theta_i, \in  \mathbb{R}/ \mathbb{Z} \}, L_2^{\infty} :=\{ (\theta_1, \theta_2 ,\frac{1}{2}) \, | \, \theta_i  \in \mathbb{R}/ \mathbb{Z}  \}.$$
We have the decomposition of $M_0 := T^3 \setminus (L_1^{\infty} \cup L_2^{\infty}) $  into connected components 
$$ \begin{array}{l} M_0 = X_1 \sqcup X_2, \\
X_1 :=\{ (\theta_1, \theta_2, \theta_3) \, | \, \theta_1, \theta_2 \in \mathbb{R}/\mathbb{Z}, 0<\theta_3<\frac{1}{2}\}, \\
X_2 :=\{(\theta_1, \theta_2, \theta_3) \, | \, \theta_1, \theta_2 \in \mathbb{R}/\mathbb{Z}, \frac{1}{2}<\theta_3<1\}.
\end{array}$$
Here any leaf of ${\cal F}|_{X_a}$ is dense in $X_a$ for each $a=1,2$.  so $\mathfrak{S}_1$ is of type III-3. Let $\omega_{\mathfrak{S}_1}$ be the smooth $1$-form defined by
$$\omega_{\mathfrak{S}_1} := {\rm cosec}(2\pi \theta_3) \omega_0 = d\theta_1+\rho d\theta_2+ {\rm cosec}(2\pi \theta_3) d\theta_3.$$
Then we see 
$$ {\rm Ker}(\omega_{\mathfrak{S}_1}) = {\rm Ker}(\omega_0), \;\omega_{\mathfrak{S}}(\dot{\phi}^t) = 1 $$  and so  $\omega_{\mathfrak{S}_1}$ is indeed the canonical $1$-form of $\mathfrak{S}_1$.
Let $\gamma_a^1, \gamma_a^2$ be closed curves in $X_a$ defined by 
$$
\gamma_a^1 :=\{(t,0,\frac{2a-1}{4}) \mid t\in \mathbb{R}/\mathbb{Z} \},
\quad 
\gamma_a^2 :=\{(0,t,\frac{2a-1}{4}) \mid t\in \mathbb{R}/\mathbb{Z} \}.
$$
Then we see that $H_1(X_a,\mathbb{Z})$ is generated by the homology classes of $\gamma_a^1, \gamma_a^2$ and that
$[\omega_{\mathfrak{S}}]([\gamma_a^1]=1,  [\omega_{\mathfrak{S}}]([\gamma_a^2])=\rho.$
Hence the period group $\Lambda_{\mathfrak{S}}$ is $\mathbb{Z} + \rho\mathbb{Z}$. For this example $\mathfrak{S}_1$, any orbits in $M_0$ is closed and so ${\cal P}_{\mathfrak{S}}$ is uncountable. However, as in Example 3.II.1, we can change $\phi_1$ to obtain an FDS$^3$ having countably infinitely many closed orbits in $M_0$ as follows. Let $V_2 := (-\rho, 1, 0), 
V_3 := (-1, 0, \sin(2\pi\theta_3))$ be vector fields on $T^3$ and we consider the smooth dynamical system $\phi_2$ defined by 
$$ \begin{array}{ll} \displaystyle{\frac{d}{dt}\phi_2^t(p)} &  := 
(1, 0, 0) + \sin(2\pi\theta_2) V_2 + \cos(2\pi\theta_3) V_3 \\
& =(1-\rho\sin(2\pi\theta_2)-\cos(2\pi\theta_3), \sin(2\pi\theta_2), \frac{1}{2}
\sin(4\pi\theta_3)).
\end{array}
$$
Then we obtain an FDS$^3$ $\mathfrak{S}_2 := (T^3, {\cal F}, \phi_2)$ of type III-3. Since $V_2$ and $V_3$ are tangent to ${\cal F}$ and so $(\dot{\phi}_1^t - \dot{\phi}_2^t)_p \in T_p{\cal F}$ for $p \in T^3$, we have $\omega_{\mathfrak{S}_2} = \omega_{\mathfrak{S}_1}, \Lambda_{\mathfrak{S}_2} = \Lambda_{\mathfrak{S}_1}$ by Lemma 1.12. We see that ${\cal P}_{\mathfrak{S}_2}$ consists of the following four closed orbits in $M_0$
$$ \begin{array}{ll}
\gamma_1 = \{ (\theta_1, 0, \frac{1}{4}) \, | \, \theta_1 \in S^1\},  & \gamma_2 = \{ (\theta_1, 0, \frac{3}{4}) \, | \, \theta_1 \in S^1\}, \\
\gamma_3 = \{ (\theta_1, \frac{1}{2}, \frac{1}{4}) \, | \, \theta_1 \in S^1 \},  &  \gamma_4 = \{ (\theta_1, \frac{1}{2}, \frac{3}{4} ) \, | \, \theta_1 \in S^1 \},
\end{array}
$$
and that $\phi_2$ is of contracting type around  $\gamma_3$ and $\gamma_4$. We replace a tubular neighborhood of $\gamma_3$ or $\gamma_4$ by the suspension of the horseshoe map. By Lemma 3.2,  we obtain an FDS$^3$ $\mathfrak{S} = (M, {\cal F}, \phi)$ of type III-3 such that  ${\cal P}_{\mathfrak{S}}$ is countably infinite, and the canonical form $\omega_{\mathfrak{S}}$ and the period group $\Lambda_{\mathfrak{S}}$ are same as $\omega_{\mathfrak{S}_i}$ and $\Lambda_{\mathfrak{S}_i}$ ($i=1,2$), respectively.\\

\begin{center}
{\bf 4.  Smooth Deligne cohomology and integration theory for FDS$^3$'s }
\end{center}

In this section, we recall the theory of smooth Deligne cohomology for an FDS$^3$ and the integration theory of Deligne cohomology classes. For general materials on smooth Deligne cohomology, we consult [Br]. For the integration theory of Deligne cohomology classes, we refer to [Ga], [GT] and [Te]. \\ 
\\
{\bf 4.1. Smooth Deligne cohomology.}  Let $\mathfrak{S} = (M, {\cal F}, \phi)$ be an FDS$^3$.  Let $X$ be a submanifold of $M_0$ obtained by removing  some finitely many closed orbits. 

Let ${\cal A}^i$ denote the sheaf of $\mathbb{C}$-valued smooth $i$-forms on $X$.  For example, an element of ${\cal A}^1(U)$ is given  in terms of a foliated local coordinate $(z,x) \in U$ by
 $$ f_1(z,x)dz + f_2(z,x)d\overline{z} + f_3(z,x)dx$$
where $f_i(z,x)$'s are $\mathbb{C}$-valued smooth functions on $U$. 

Let $\Lambda$ be a subgroup of the additive group $\mathbb{R}$. For a non-negative integer $n$, we set $\Lambda(n) := (2\pi \sqrt{-1})^n \Lambda$.\\
\\
{\bf Definition 4.1.1.} Let $n$ be an integer with $1\leq n \leq 3$. We define the 
{\em smooth Deligne complex}  $\Lambda(n)_{\mathscr{D}}$  on $X$ by 
$$ \Lambda(n)_{\mathscr{D}} :  \; \Lambda(n) \rightarrow {\cal A}^0 \stackrel{d}{\rightarrow} {\cal A}^1 \stackrel{d}{\rightarrow} \cdots \stackrel{d}{\rightarrow}
{\cal A}^{n-1},$$
where $\Lambda(n)$ is put in degree $0$ and $d$ denotes the differential. 
For an integer $q \geq 0$, the $q$-th {\em smooth Deligne cohomology group} with coefficients in $\Lambda(n)_{\mathscr{D}}$ is  defined to be the  
$q$-th hypercohomology group of the complex $\Lambda(n)_{\mathscr{D}}$, denoted by $H^q_{\mathscr{D}}(M; \Lambda(n))$:
$$ H^q_{\mathscr{D}}(X; \Lambda(n)):= \mathbb{H}^q(X; \Lambda(n)_{\mathscr{D}}). $$
In particular, when $\Lambda$ is the period group $\Lambda_{\mathfrak{S}}$, we call $H^q_{\mathscr{D}}(X; \Lambda_{\mathfrak{S}}(n))$ the {\em FDS$^3$-Deligne cohomology groups} of $\mathfrak{S}$.\\
\\
We compute the smooth Deligne cohomology groups as  \v{C}ech hypercohomology groups of an open covering ${\cal U} = \{U_a\}_{a \in I}$ of $X$ with coefficients in $\Lambda(n)_{\mathscr{D}}$
$$ H^q_{\mathscr{D}}(M; \Lambda(n)) = \mathbb{H}^q({\cal U}; \Lambda(n)_{\mathscr{D}}),$$
where the open covering ${\cal U}$ is taken so that all non-empty intersections $U_{a_0\cdots a_j} := U_{a_0}\cap \dots \cap U_{a_j}$ are contractible. So a \v{C}ech cocycle representing an element of $H^n_{\mathscr{D}}(M; \Lambda(n))$ is of the form
$$ (\lambda_{a_0 \dots a_n}, \theta^0_{a_0 \dots a_{n-1}}, \dots , 
\theta^{n-1}_{a_0}) \in C^n(\Lambda(n)) \oplus C^{n-1}({\cal A}^0)\oplus \cdots \oplus C^0({\cal A}^{n-1}) $$
which satisfies the cocycle condition
$$ \delta(\theta^0_{a_0 \dots a_{n-1}}) + (-1)^n \lambda_{a_0 \dots a_n} = 0, \;\; \delta(\theta^i_{a_0 \dots a_{n-1-i}}) +(-1)^{n-i} d\theta^{i-1}_{a_0 \dots a_{n-i}} = 0 \; (i \geq 1),$$
where $\delta$ is the \v{C}ech differential with respect to the open covering ${\cal U}$.\\
\\
{\bf Example 4.1.2.} Let $f$ be an FDS$^3$-meromorphic function on $\mathfrak{S}$ whose zeros and poles are lying along closed orbits $\gamma_1, \dots , \gamma_N$. Let $X := M_0$. Let $\log_a f$ denote a branch of $\log f$ on $U_a$. Then the \v{C}ech cocycle
$$ (n_{a_0a_1}, \log_{a_0} f), \; n_{a_0a_1} = (\delta \log f)_{a_0a_1} \in \mathbb{Z}(1) $$
determines the cohomology class of $H^1_{\mathscr{D}}(X; \mathbb{Z}(1))$, by which we denote $c(f)$. 
\\
\\
{\bf Example 4.1.3.} Let $\omega_{\mathfrak{S}}$ be the canonical $1$-form of $\mathfrak{S}$ (cf. Definition 1.11). Fix a base point $p_0 \in X$. For each $U_a (a \in I)$, we choose a point $p_a \in U_a$ and a path $\gamma_a$ in $U_a$ from $p_0$ to $p_a$. For $p \in U_a$, we set
$$ f_{\omega_{\mathfrak{S}},a}(p) := 2 \pi \sqrt{-1} \int_{\gamma_p \cdot \gamma_a} \omega_{\mathfrak{S}},$$
where $\gamma_p$ is a path from $p_a$ to $p$ inside $U_a$. Since $U_a$ is contractible and $\omega_{\mathfrak{S}}$ is closed, $f_{\omega_{\mathfrak{S},a}}$ is a smooth function on $U_a$ which is independent of the choice of $\gamma_p$. Then the \v{C}ech cocycle 
$$ (\lambda_{a_0a_1}, f_{\omega_{\mathfrak{S}},a_0}), \; \lambda_{a_0a_1} = (\delta f_{\omega_{\mathfrak{S}}})_{a_0a_1} \in \Lambda_{\mathfrak{S}}(1)$$ 
defines the cohomology class of $H^1_{\mathscr{D}}(X; \Lambda_{\mathfrak{S}}(1))$, by which we denote $c(\omega_{\cal S})$. 
The class $c(\omega_{\cal S})$ is independent of the choices of $ p_a$, $\gamma_a$ and $\gamma_p$. In fact, let $p_a'$, $\gamma_a'$ and $\gamma_p'$ be different choices of a point in $U_a$, a path from $p_0$ to $p_a'$ and a path from $p_a'$ to $p$, respectively, and let $f_{\omega_{\mathfrak{S}}}'$ and $\lambda_{a_0a_1}'$ be defined as above using $p_a', \gamma_a'$ and $\gamma_p'$. Since $d f_{\omega_{\mathfrak{S},a}}' = d f_{\omega_{\mathfrak{S},a}} = \omega_{\mathfrak{S}}$, $f_{\omega_{\mathfrak{S},a}}' - f_{\omega_{\mathfrak{S}},a} = \lambda_a$ is a constant on $U_a$. Then we have 
$$\lambda_a = f_{\omega_{\mathfrak{S},a}}'(p) - f_{\omega_{\mathfrak{S},a}}(p) = 2\pi\sqrt{-1} \int_{(\gamma_p \cdot \gamma_a)^{-1} \cdot (\gamma_p' \cdot \gamma_a')} \omega_{\mathfrak{S}} \in \Lambda_{\mathfrak{S}}(1)$$
and 
$$\begin{array}{ll} \lambda_{a_0a_1}' - \lambda_{a_0a_1} & = (\delta f_{\omega_{\mathfrak{S}}}')_{a_0a_1} - (\delta f_{\omega_{\mathfrak{S}}}')_{a_0a_1}\\
 &  = f_{\omega_{\mathfrak{S}, a_1}}'- f_{\omega_{\mathfrak{S}, a_0}}' -( f_{\omega_{\mathfrak{S}, a_1}} - f_{\omega_{\mathfrak{S}, a_0}})\\
& = (f_{\omega_{\mathfrak{S}, a_1}}' - f_{\omega_{\mathfrak{S}, a_1}}) - (f_{\omega_{\mathfrak{S}, a_0}}' - f_{\omega_{\mathfrak{S}, a_0}}) \\
& = \lambda_{a_1} - \lambda_{a_0} \\
& = \delta(\lambda)_{a_0a_1}
\end{array}$$
and hence $c(\omega_{\cal S})$ is independent of the choices of $p_a$, $\gamma_a$ and $\gamma_p$, but, it depends on the choice of $p_0$.\\
\\
{\bf Definition 4.1.4.} For $1\leq n \leq 3$, we define the {\em $n$-curvature homomorphim}
$$ \Omega : H^n_{\mathscr{D}}(X; \Lambda(n)) \longrightarrow {\cal A}^n(X)$$
by
$$ \Omega(c)|_{U_a} := d\theta_a^{n-1} $$
 for $c = [(\lambda_{a_0 \dots a_n}, \dots , \theta^{n-1}_{a_0})]$. 
\\

When $\Lambda$ is a subring of $\mathbb{R}$, the smooth Deligne cohomology groups are equipped with the cup product, which is induced by the product  on the smooth Deligne complexes
$$ \Lambda(n)_{\mathscr{D}} \otimes \Lambda(n')_{ \mathscr{D}}\longrightarrow \Lambda(n+n')_{\mathscr{D}}  \leqno{(4.1.5)}$$
defined by
$$ x \cup y = \left\{ \begin{array}{ll}
xy & \; {\rm deg}(x) = 0,\\
x \wedge dy & \; {\rm deg}(x) > 0 \; {\rm and}\; {\rm deg}(y) = n',\\
0 & \; \mbox{otherwise}.
\end{array}\right.
 \leqno{(4.1.6)}
$$
For our purpose, we extend the  product (4.1.5) for the case where $\Lambda$ is a subring of $\mathbb{R}$ and $\Lambda'$ is a $\Lambda$-submodule of $\mathbb{R}$ as follows. Namely, by 
the same formula as in (4.1.6), we have the product
$$\Lambda(n)_{\mathscr{D}} \otimes \Lambda'(n')_{ \mathscr{D}}\longrightarrow \Lambda'(n+n')_{\mathscr{D}}  $$
which induces the cup product on the smooth Deligne cohomology groups
$$ H^n_{\mathscr{D}}(X;\Lambda(n)) \otimes H^{n'}_{\mathscr{D}}(X;\Lambda'(n')) \longrightarrow H^{n+n'}_{\mathscr{D}}(X;\Lambda'(n+n')). \leqno{(4.1.7)}$$
\\
{\bf 4.2. Integration theory.} As in the subsection 4.1, let $\mathfrak{S} = (M, {\cal F},\phi)$ be an FDS$^3$ and let $X$ be a submanifold of $M_0$ obtained by removing  some  finitely many closed orbits.
Let $n$ be an integer with $1\leq n \leq 3$, and let $c \in H^n_{\mathscr{D}}(X;\Lambda(n))$. Let $Y$ be an $(n-1)$-dimesnional closed submanifold of $X$. We shall define  a paring 
$\int_Y c$, which takes values in $\mathbb{C}$ mod $\Lambda(n)$, as follows.

First,  we fix an open covering ${\cal U}=\{ U_{a} \}_{a \in I}$ of $X$ such that all non-empty intersections $U_{a_0\cdots a_j} := U_{a_0}\cap \dots \cap U_{a_j}$ are contractible and choose  a \v{C}ech representative 
cocycle $(\lambda_{a_0 \dots a_n}, \theta^0_{a_0 \dots a_{n-1}}, \dots , 
\theta^{n-1}_{a_0})$ of $c$. 
Second, we choose a smooth finite  triangulation $K = \{ \sigma \}$ of $Y$ and an index map $ \iota : K \rightarrow I$ satisfying $\sigma \subset U_{\iota(\sigma)}$. For $i = 0,\dots, n-1$, we define the set $F_K(i)$ of flags of simplices
$$ F_K(i) := \{ \vec{\sigma} = (\sigma^{n-1-i},\dots , \sigma^{n-1}) \, | \, \sigma^j \in K, \dim \sigma^j = j, \sigma^{n-1-i} \subset \cdots \subset  \sigma^{n-1}\}. \leqno{(4.2.1)}$$
Then, we define  the {\em integral} of $c$ over $Y$ by 
$$
\int_Y c:=\sum_{i=0}^{n-1} \sum_{\vec{\sigma} \in F_K(i)} \int_{\sigma^{n-1-i}} 
\theta^{n-1-i}_{\iota(\sigma^{n-1}) \iota(\sigma^{n-2}) \dots \iota(\sigma^{n-1-i})} 
\mod \Lambda(n),
\leqno{(4.2.2)}
$$
which is  proved to be  independent of all choices ([GT; Theorem 3.4 (i)]). \\
\\
The following Stokes-type formula was shown by Gawedzki ([Ga]) when $Y$ is $2$-dimensional. We refer to [GT] and [Te] for more general statements and proofs. \\
\\
{\bf Theorem 4.2.3} (cf. [Ga], [GT; Theorem 3.4 (ii)], [Te; Proposition 5.5]). {\em If there is an $n$-dimensional submanifold $Z$ of $X$ whose boundary is $\partial Z=Y$, we have}
$$ \int_Y c = \int_Z \Omega(c) \;\; \mbox{mod} \; \Lambda(n). $$
\\

\begin{center}
{\bf 5. Local symbols and Hilbert type reciprocity law}
\end{center}

In this section, we introduce a local symbol by using the integral  of a certain FDS-Deligne cohomology class along a torus and show the Hilbert type reciprocity law.\\
\\
{\bf 5.1. Local symbols.} Let $\mathfrak{S} = (M, {\cal F},\phi)$ be an FDS$^3$. Let $f$ and $g$ be FDS$^3$-meromorphic functions on $\mathfrak{S}$ whose zeros and poles lie along $\gamma_1, \dots , \gamma_N \in {\cal P}_{\mathfrak{S}}$. We set $X := M_0 \setminus  \cup_{i=1}^N \gamma_i$. For $\gamma \in \overline{{\cal P}_{\mathfrak{S}}}$, let $V(\gamma)$ denote a tubular neighborhood of $\gamma$ and we denote by $T(\gamma)$ the boundary of $V(\gamma)$.

 As in Examples 4.1.2, we have the smooth Deligne cohomology classes
$$ c(f) = [(m_{a_0a_1}, \log_{a_0} f)], \, c(g) = [(n_{a_0a_1}, \log_{a_0} g)] \in H^1_{\mathscr{D}}(X; \mathbb{Z}(1))$$
and, as in Example 4.1.3, we have the FDS$^3$-Deligne cohomology class
$$ c(\omega_{\mathfrak{S}}) = [(\lambda_{a_0a_1}, f_{\omega_{\cal S},a_0})] \in H^1_{\mathscr{D}}(X; \Lambda_{\mathfrak{S}}(1)).$$
By the product in (4.1.7) applied to the case that $\Lambda = \mathbb{Z}$ and $\Lambda' = \Lambda_{\mathfrak{S}}$, we have the $3$rd FDS$^3$-Deligne cohomology class
$$ c(f) \cup c(g) \cup c(\omega_{\mathfrak{S}}) \in H^3_{\mathscr{D}}(X;\Lambda_{\mathfrak{S}}(3)).$$
\\
{\bf Definition 5.1.1.} We define the {\em local symbol} $\langle f, g \rangle_{\gamma}$ of $f, g$ along $\gamma$ by
$$ \langle f, g \rangle_{\gamma} := \int_{T(\gamma)} c(f) \cup c(g) \cup c(\omega_{\mathfrak{S}})  \;\; \mbox{mod} \; \Lambda_{\mathfrak{S}}(3).$$
We note that the integral of the r.h.s. is finite since $T(\gamma)$ is compact and that it is independent of a choice of $V(\gamma)$ by the Stokes theorem.\\
\\
{\bf Theorem 5.1.2.} {\em Notations being as above, the FDS$^3$-Deligne cohomology class $c(f)\cup c(g) \cup c(\omega_{\mathfrak{S}})$ is represented by the \v{C}ech cocycle}
$$ (m_{a_0a_1}n_{a_1a_2}\lambda_{a_2a_3}, \, m_{a_0a_1}n_{a_1a_2} f_{\omega_{\mathfrak{S}}, a_0a_1a_2}, \, m_{a_0a_1}\log_{a_1}g \, \omega_{\mathfrak{S}}, \, \log_{a_0} f d\log g \wedge \omega_{\mathfrak{S}}). $$
{\em For $\gamma \in {\cal P}_{\mathfrak{S}} $, the local symbol $\langle f, g \rangle_{\gamma}$ is given by}
$$ \langle f, g \rangle_{\gamma} = \int_{T(\gamma)} \log(f) d\log(g) \wedge \omega_{\mathfrak{S}} + \int_{\mathfrak{m}} d\log(f)  \int_{\mathfrak{l}} \log(g) \omega_{\mathfrak{S}} \; \; \mbox{mod}\; \Lambda_{\mathfrak{S}}(3),$$
{\em where $\mathfrak{m}$ and $\mathfrak{l}$ denote a meridian and longitude on $T(\gamma)$, respectively.}\\
\\
{\em Proof.} The first assertion follows from the definition (4.1.6) of the cup product on Deligne cohomology groups. We set for simplicity
$$ \left\{
\begin{array}{ll}
\theta^0_{a_0a_1a_2} & = m_{a_0a_1}n_{a_1a_2} f_{\omega_{\mathfrak{S}}, a_0a_1a_2}, \\
\theta^1_{a_0a_1} & = m_{a_0a_1}\log_{a_1}g\, \omega_{\mathfrak{S}}, \\
\theta^2_{a_0} & = \log_{a_0} f d\log g \wedge \omega_{\mathfrak{S}}.
\end{array} \right. 
$$
To prove the second formula for the local symbol, we make suitable choices of a triangulation $K$ of $T(\gamma)$, an open covering ${\cal U}$, representatives of  Deligne cohomology classes $c(f), c(g), c(\omega_{\mathfrak{S}})$, and the index map $\iota$. First, let $K = \{ \sigma \}$ be a triangulation of the $2$-dimensional torus $T(\gamma)$ by eighteen triangles as in the following figure:

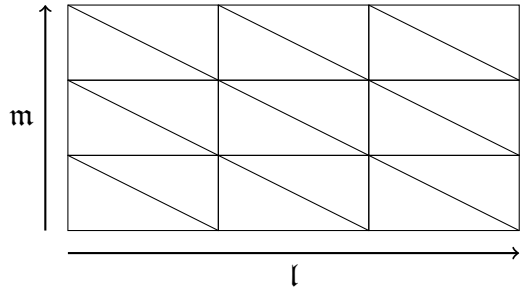
\begin{figure}[H]
\centering
\begin{tikzpicture}
\draw (0,2) rectangle (2,3);\draw (2,2) rectangle (4,3);\draw (4,2) rectangle (6,3);
\draw (0,1) rectangle (2,2);\draw (2,1) rectangle (4,2);\draw (4,1) rectangle (6,2);
\draw (0,0) rectangle (2,1);\draw (2,0) rectangle (4,1);\draw (4,0) rectangle (6,1);
\draw (0,1) -- (2,0);
\draw (0,2) -- (4,0);
\draw (0,3) -- (6,0);
\draw (2,3) -- (6,1);
\draw (4,3) -- (6,2);
\draw[thick,->] (0,-0.3) -- (6,-0.3) node[pos=0.5, below] {$\mathfrak{l}$};
\draw[thick,->] (-0.3,0) -- (-0.3,3) node[pos=0.5, left] {$\mathfrak{m}$};
\end{tikzpicture}
\end{figure}
We choose an open covering ${\cal U} = \{ U_a \}_{a \in I}$ of $X$ such that nine $U_a$'s in ${\cal U}$ cover $T(\gamma)$ and are indexed,  as in the following figure

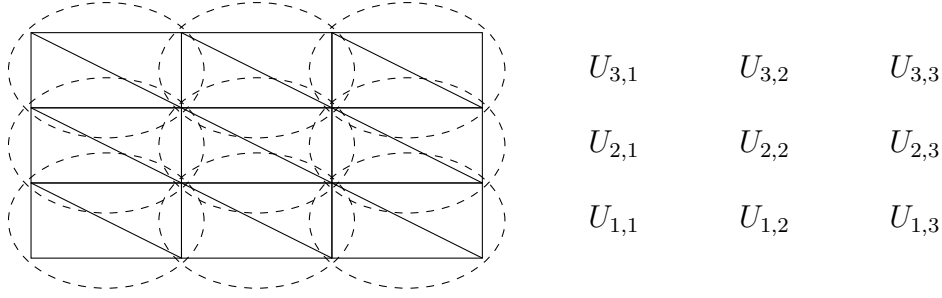
\begin{figure}[H]
\centering
\begin{tikzpicture}
\draw (0,2) rectangle (2,3);\draw (2,2) rectangle (4,3);\draw (4,2) rectangle (6,3);
\draw (0,1) rectangle (2,2);\draw (2,1) rectangle (4,2);\draw (4,1) rectangle (6,2);
\draw (0,0) rectangle (2,1);\draw (2,0) rectangle (4,1);\draw (4,0) rectangle (6,1);
\draw (0,1) -- (2,0);
\draw (0,2) -- (4,0);
\draw (0,3) -- (6,0);
\draw (2,3) -- (6,1);
\draw (4,3) -- (6,2);
\draw[dashed] (1,0.5) ellipse (1.3cm and 0.9cm);
\draw[dashed] (3,0.5) ellipse (1.3cm and 0.9cm);
\draw[dashed] (5,0.5) ellipse (1.3cm and 0.9cm);
\draw[dashed] (1,1.5) ellipse (1.3cm and 0.9cm);
\draw[dashed] (3,1.5) ellipse (1.3cm and 0.9cm);
\draw[dashed] (5,1.5) ellipse (1.3cm and 0.9cm);
\draw[dashed] (1,2.5) ellipse (1.3cm and 0.9cm);
\draw[dashed] (3,2.5) ellipse (1.3cm and 0.9cm);
\draw[dashed] (5,2.5) ellipse (1.3cm and 0.9cm);
\end{tikzpicture}
\begin{tikzpicture}
\draw[white] (1,0.5) ellipse (1.3cm and 0.9cm);
\draw[white] (3,0.5) ellipse (1.3cm and 0.9cm);
\draw[white] (5,0.5) ellipse (1.3cm and 0.9cm);
\draw[white] (1,1.5) ellipse (1.3cm and 0.9cm);
\draw[white] (3,1.5) ellipse (1.3cm and 0.9cm);
\draw[white] (5,1.5) ellipse (1.3cm and 0.9cm);
\draw[white] (1,2.5) ellipse (1.3cm and 0.9cm);
\draw[white] (3,2.5) ellipse (1.3cm and 0.9cm);
\draw[white] (5,2.5) ellipse (1.3cm and 0.9cm);
\node at (1,0.5) {$U_{1,1}$};
\node at (3,0.5) {$U_{1,2}$};
\node at (5,0.5) {$U_{1,3}$};
\node at (1,1.5) {$U_{2,1}$};
\node at (3,1.5) {$U_{2,2}$};
\node at (5,1.5) {$U_{2,3}$};
\node at (1,2.5) {$U_{3,1}$};
\node at (3,2.5) {$U_{3,2}$};
\node at (5,2.5) {$U_{3,3}$};
\end{tikzpicture}
\end{figure}

We define the index map $\iota : K \rightarrow I$ in the manner that one vertex, three edges and two triangles in $U_a$ are sent to $a \in I$ as in the following figure:

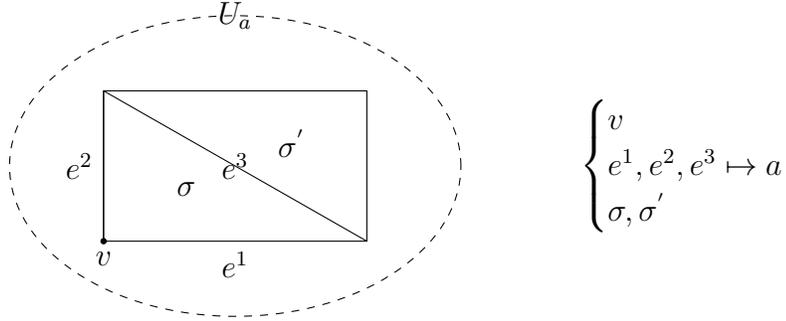
\begin{figure}[H]
\centering
\begin{tikzpicture}
\draw 
(0,0)  -- node[anchor=south, below] {$e^{1}$}
(3.5,0) -- (3.5,2) -- (0,2) -- cycle; 
\draw (0,2) -- (3.5,0) node[pos=0.5] {$e^{3}$};
\draw[dashed] (1.75,1) ellipse (3cm and 2cm);
\filldraw 
(0,2) -- node[pos=0.5, left] {$e^{2}$}
(0,0) circle (1pt) node[align=left, below] {$v$}
(1.75,3) circle (0pt) node[align=center] {$U_{a}$}
(1.1,0.7) circle (0pt) node[align=center] {$\sigma$}
(2.5,1.3) circle (0pt) node[align=center] {$\sigma^{'}$};
\end{tikzpicture}
\begin{tikzpicture}
\draw[white] (1.75,1) ellipse (3cm and 2cm);
\node at (1.75,1) {$
    \begin{cases}
    v\\
    e^{1}, e^{2}, e^{3} \mapsto  a\\
    \sigma,\sigma^{'}
    \end{cases}
    $};
\end{tikzpicture}
\end{figure}

We choose the representatives of $c(f), c(g)$ and $c(\omega_{\mathfrak{S}})$ as follows:
$$ \begin{array}{l}
c(f) = [(m_{a_0a_1}, \log_{a_0} f)], \; m_{a_0a_1} = (\delta \log f)_{a_0a_1}\\
m_{a_0a_1} = \left\{ \begin{array}{ll} \displaystyle{\int_{\mathfrak{m}} d\log f}  & \; \cdots \; a_0 = (3,i), a_1 = (1,j)\\
 \displaystyle{- \int_{\mathfrak{m}} d\log f } & \; \cdots \; a_0 = (1,i), a_1 = (3,j)\\
0  & \; \cdots \;  \mbox{otherwise}, \end{array} \right.\\
c(g) = [(n_{a_0a_1}, \log_{a_0} g)], \; n_{a_0a_1} = (\delta \log g)_{a_0a_1} \\
n_{a_0a_1} = \left\{ \begin{array}{ll} \displaystyle{\int_{\mathfrak{m}} d\log g}  & \; \cdots \; a_0 = (3,i), a_1 = (1,j)\\
 \displaystyle{- \int_{\mathfrak{m}} d\log g } & \; \cdots \; a_0 = (1,i), a_1 = (3,j)\\
0  & \; \cdots \;  \mbox{otherwise}, \end{array} \right.\\
c(\omega_{\mathfrak{S}}) = [(\lambda_{a_0a_1}, f_{a_0}^{\omega_{\mathfrak{S}}})], \; \lambda_{a_0a_1} = (\delta f_{\omega_{\mathfrak{S}}})_{a_0a_1} \\
\lambda_{a_0a_1} = \left\{ \begin{array}{ll} \displaystyle{\int_{\mathfrak{l}} \omega_{\mathfrak{S}}}  & \; \cdots \; a_0 = (i,3), a_1 = (j,1)\\
 \displaystyle{- \int_{\mathfrak{l}}  \omega_{\mathfrak{S}}} & \; \cdots \; a_0 = (i,1), a_1 = (j,3)\\
0  & \; \cdots \;  \mbox{otherwise}, \end{array} \right.
\end{array}
$$
where $i, j \in \{ 1,2,3 \}$.  

Let $F_K(i)$ be the set of flags of simplices for $i = 0,1,2$ as in (4.2.1). By (4.2.2) and Definition 5.1.1, we have
$$ \begin{array}{ll} 
\langle f, g \rangle_{\gamma}  &   =    \displaystyle{ \sum_{\vec{\sigma} \in F_K(0)} \int_{\sigma^2} \theta^2_{\iota(\sigma^2)} } \\
                              & \; + \displaystyle{ \sum_{\vec{\sigma} \in F_K(1)} \int_{\sigma^1}  \theta^1_{\iota(\sigma^2) \iota(\sigma^1)}  }  \\
                         & \; +  \displaystyle{ \sum_{\vec{\sigma} \in F_K(2)} \int_{\sigma^0} \theta^0_{\iota(\sigma^2)\iota(\sigma^1)\iota(\sigma^0)}  }.
   \end{array}                       \leqno{(5.1.2.1)}
$$
For the 1st term of the r.h.s. of (5.1.2.1), we have
$$ \begin{array}{ll} \displaystyle{ \sum_{\vec{\sigma} \in F_K(0)} \int_{\sigma^2} \theta^2_{\iota(\sigma^2)} } &   = \displaystyle{ \sum_{a \in I} \sum_{\sigma^2 \subset U_a} \int_{\sigma^2} \log_{a} f d\log g \wedge \omega_{\mathfrak{S}} }\\
  & = \displaystyle{ \int_{T(\gamma)} \log f d\log g \wedge \omega_{\mathfrak{S}}}.
  \end{array}  \leqno{(5.1.2.2)}
  $$
 For the 2nd term of the r.h.s. of (5.1.2.1), we note that $m_{\iota(\sigma^2) \iota(\sigma^1)} = 0$ unless $(\iota(\sigma^2), \iota(\sigma^1)) = ((3,i), (1,i))$ for $i = 1,2,3$, by our choices of $m_{ab}$'s and $\iota$. Hence we have
 $$ \begin{array}{ll}
\displaystyle{ \sum_{\vec{\sigma} \in F_K(1)} \int_{\sigma^1}  \theta^1_{\iota(\sigma^2) \iota(\sigma^1)} } 
& =  \displaystyle{  \sum_{a \in I} \sum_{\sigma^1 \subset \sigma^2 \subset U_a} \int_{\sigma^1} m_{\iota(\sigma^2) \iota(\sigma^1)} \log_{\iota(\sigma^1)} g \, \omega_{\mathfrak{S}} }\\
& = \displaystyle{ \sum_{i=1}^3 \sum_{\iota(\sigma^1) \subset U_{3,i}} \int_{\sigma^1} m_{(3,i)(1,i)} \log_{\iota(\sigma^1)} g \, \omega_{\mathfrak{S}} }\\
& = \displaystyle{ \int_{\mathfrak{m}} d\log f \int_{\mathfrak{l}} \log g \, \omega_{\mathfrak{S}}. }
\end{array}
\leqno{(5.1.2.3)} $$
For the 3rd term of the r.h.s. of (5.1.2.1), we note that  $m_{\iota(\sigma^2)\iota(\sigma^1)}n_{\iota(\sigma^1)\iota(\sigma^0)} = 0$ for any flags,  by our choices of $m_{ab}$'s, $n_{ab}$'s and $\iota$. Hence we have
$$\displaystyle{ \sum_{\vec{\sigma} \in F_K(2)} \int_{\sigma^0}  \theta^0_{\iota(\sigma^2) \iota(\sigma^1) \iota(\sigma^0)}  = 0. } \leqno{(5.1.2.4)} $$
 Getting $(5.1.2.2) \sim (5.1.2.4)$ together, we obtain the second formula for $\langle f, g \rangle_{\gamma}$. $\;\; \Box$
 \\
 \\
 {\bf Remark 5.1.3.} (1) Brylinski and McLaughlin showed a formula for the holonomy of a certain gerbe along a torus ([BM; Theorem 3.6]). We note that their formula has a form similar to ours above and that their method of the computation is different from ours.\\
(2) Bloch pointed out to us that the formula in Theorem 5.1.2 has a form similar to a formula for
the regulator of an elliptic curve (cf. [Be; 4.2], [DW; 1.10])).
\\
\\
{\bf 5.2. Hilbert type reciprocity law.}  We keep the same notations as in the subsection 5.1. We show a geometric analogue for our local symbol of the reciprocity law for the Hilbert symbol in a global field. \\
\\
{\bf Theorem 5.2.1.} {\em We have }
$$
\sum_{\gamma \in \overline{{\cal P}_{\mathfrak{S}}}} \langle f,g \rangle_{\gamma}=0 \mod \Lambda_{\mathfrak{S}}(3).
$$
\\
{\em Proof.} We set $Z := X \setminus (\cup_{i=1}^r V(\gamma_{i}^{\infty}) \cup \cup_{i=1}^N V(\gamma_i))$. Noting $\partial Z = \cup_{i=1}^r T(\gamma_{i}^{\infty}) \cup \cup_{i=1}^N T(\gamma_i)$,   we have, by Theorem 4.2.3,
$$ \begin{array}{ll}\displaystyle{\sum_{\gamma \in \overline{{\cal P}_{\mathfrak{S}}}} \langle f,g \rangle_{\gamma}} & = \displaystyle{\sum_{i=1}^r  \langle f,g \rangle_{\gamma_i^{\infty}} + \sum_{i=1}^N  \langle f,g \rangle_{\gamma_i} }\\
  & = \displaystyle{ \sum_{i=1}^r \int_{T(\gamma_i^{\infty})} c(f) \cup c(g) \cup c(\omega_{\mathfrak{S}}) }\\
 &   \;\;\;\;\; \displaystyle{+ \sum_{i=1}^N \int_{T(\gamma_i)} c(f) \cup c(g) \cup c(\omega_{\mathfrak{S}}) }
\;\; \mbox{mod}\; \Lambda_{\mathfrak{S}}(3)\\
 & = \displaystyle{ \int_{\partial Z} c(f) \cup c(g) \cup c(\omega_{\mathfrak{S}}) } \;\; \mbox{mod}\; \Lambda_{\mathfrak{S}}(3) \\
 & = \displaystyle{ \int_{Z} \Omega(c(f) \cup c(g) \cup c(\omega_{\mathfrak{S}}))  \;\; \mbox{mod}\; \Lambda_{\mathfrak{S}}(3).}
 \end{array}
$$
By Definition 4.1.4 and Theorem 5.1.2, we have
$$ \Omega(c(f) \cup c(g) \cup c(\omega_{\mathfrak{S}}))|_{U_a} = d(\log_{a} f d\log g \wedge \omega_{\mathfrak{S}}).$$
Taking a foliated local coordinate $(z, x)$ on $U_a$, we have
$$ \begin{array}{ll}
d(\log_{a} f d\log g \wedge \omega_{\mathfrak{S}}) & = d(\log_{a} f d\log g) \wedge  \omega_{\mathfrak{S}} \;\; (\mbox{since}\; \omega_{\mathfrak{S}}\; \mbox{is closed}) \\
                  & = d\log f \wedge  d\log g \wedge \omega_{\mathfrak{S}} \\
 & =  \displaystyle{ (\frac{f_z}{f} dz +   \frac{f_x}{f} dx) \wedge   (\frac{g_z}{g} dz +   \frac{g_x}{g} dx)   
\wedge \omega_{\mathfrak{S}}   } \\
 & \;\;\;\;\;\;\;\;\;\; \;\;\;\;\;\; (\mbox{since}\; \log f, \log g \; \mbox{are holomorphic on} \; U_a)\\
  & =\displaystyle{ \frac{f_z g_x - f_x g_z}{fg}  dz \wedge dx \wedge \omega_{\mathfrak{S}} }\\
  & = 0 \;\;  (\mbox{by}\; \omega_{\mathfrak{S}} = h(x) dx)
\end{array}
$$
and hence the assertion follows. $\;\; \Box$\\
\\
{\bf Remark 5.2.2.} (1) Our method to introduce local symbols and to show the reciprocity law  may be regarded as a generalization to an FDS$^3$ of  Deligne, Bloch and Beilinson's interpretation of the tame symbol on a Riemann surface using the holomorphic Deligne cohomology ([Be], [Bl], [Dl]). \\
(2) It would be interesting to generalize our local symbol to  multiple symbols $\langle f_1,\dots, f_n \rangle_{\gamma}$ for several FDS$^3$-meromorphic functions $f_i$'s by using the Massey products in the smooth Deligne cohomology and the iterated integrals, as the tame symbol on a Riemann surface was generalized to polysymbols in [MT]. \\
\\
{\em Ackowledgement.} We would like to thank Christopher Deninger, Spencer Bloch, Atsushi Katsuda and Jonas Stelzig for useful communications. J.K. is supported by Grants-in-Aid for JSPS Fellow (DC1) Grant Number 17J02472. M.M. is  supported by Grant-in-Aid for Scientific Research (B) Grant Number JP17H02837. Y.T. is  supported by Grants-in-Aid for Scientific Research (C) Grant Number 17K05243 and by JST CREST Grant Number  JPMJCR14D6, Japan.
\\
\begin{flushleft}
{\bf References}\\
{[AT]} E. Artin, J. Tate, Class field theory, 2nd edition. Advanced Book Classics. Addison-Wesley Publishing Company, Advanced Book Program, Redwood City, CA, 1990.  \\
{[Be]} A. Beilinson, Higher regulators and values of $L$-functions, J. Soviet Math., {\bf 30} (1985), 2036--2070.\\
{[Bl]} S. Bloch, The dilogarithm and extensions of Lie algebras,  In: Algebraic K -theory, Evanston 1980 (Proc. Conf., Northwestern Univ., Evanston, Ill., 1980),  pp. 1--23, Lecture Notes in Math., {\bf 854}, Springer, Berlin-New York, 1981.\\
{[Br]} J.-L. Brylinski, Loop spaces, characteristic classes and geometric quantization, Progress in Mathematics, {\bf 107}. Birkhauser Boston, Inc., Boston, MA, 1993. \\
{[BM]} J.-L. Brylinski, D.A. McLaughlin, The geometry of two-dimensional symbols, K-Theory  {\bf 10}  (1996),  no. 3, 215--237. \\
{[Cl]} D. Calegari, Foliations and the geometry of 3-manifolds, Oxford Mathematical Monographs. Oxford University Press, Oxford, 2007.\\
{[CC]} A. Candel, L. Conlon,  Foliations. I,  Graduate Studies in Mathematics, {\bf 23}, American Mathematical Society, Providence, RI, 2000.\\
{[CtC]} J. Cantwell, L. Conlon, Tischler fibrations of open, foliated sets, Ann. Inst. Fourier, Grenoble {\bf 31}, 2(1981), 113--135.\\
{[Dl]} P. Deligne, Le symbole mod\'{e}r\'{e}, Inst. Hautes Etudes Sci. Publ. Math.  No. {\bf 73}  (1991), 147--181.\\
{[D1]} C. Deninger, Some analogies between number theory and dynamical systems on foliated spaces, Doc. Math. J. DMV. Extra Volume ICM I (1998), 23-46.\\
{[D2]} C. Deninger, On dynamical systems and their possible significance for arithmetic geometry, In: Regulators in analysis, geometry and number theory,  29--87, 
Progr. Math., {\bf 171}, Birkhauser Boston, Boston, MA, 2000. \\
{[D3]} C. Deninger, Number theory and dynamical systems on foliated spaces, Jahresber. Deutsch. Math.-Verein.  {\bf 103}  (2001),  no. 3, 79--100.\\
{[D4]} C. Deninger, A note on arithmetic topology and dynamical systems,  In: Algebraic number theory and algebraic geometry,  99--114, 
Contemp. Math., {\bf 300}, Amer. Math. Soc., Providence, RI, 2002. \\
{[D5]} C. Deninger, Arithmetic Geometry and Analysis on Foliated Spaces, preprint.\\
{[D6]} C. Deninger, Analogies between analysis on foliated spaces and arithmetic geometry, Groups and analysis, London Math. Soc. Lecture Note Ser. {\bf 354}, 174--190, Cambridge Univ. Press 2008.\\
{[D7]} C. Deninger, A dynamical systems analogue of Lichtenbaum's conjectures on special values of Hasse-Weil zeta functions, preprint.\\
{[DW]} C. Deninger, K. Wingberg, On the Beilinson conjectures for elliptic curves with complex multiplication, In: Beilinson's conjectures on special values of $L$-functions,  249--272, Perspect. Math., 4, Academic Press, Boston, MA, 1988.\\
{[F]} M. Farber, Topology of closed one-forms, Mathematical Surveys and Monographs, {\bf 108}. Amer. Math. Soc., Providence, RI, 2004. \\
{[Ga]} K. Gawedzki, Topological actions in two-dimensional quantum field theories,  Nonperturbative quantum field theory (Cargese, 1987),  101--141, NATO Adv. Sci. Inst. Ser. B Phys., {\bf 185}, Plenum, New York, 1988. \\
{[Gh]} \'{E}. Ghys, Laminations par surfaces de Riemann,  In: Dynamique et g\'{e}om\'{e}trie complexes (Lyon, 1997),  49--95, Panor. Syntheses, {\bf 8}, Soc. Math. France, Paris, 1999.\\
{[GT]} K. Gomi, Y. Terashima, A fiber integration formula for the smooth Deligne cohomology, Internat. Math. Res. Notices  2000,  no. {\bf 13}, 699--708.\\
{[He]} J. Hempel, $3$-manifolds, 3 -Manifolds,  Ann. of Math. Studies, No. {\bf 86}. Princeton University Press, Princeton, N. J.; University of Tokyo Press, Tokyo, 1976. \\
{[HM]} T. Horiuchi, Y. Mitsumatsu, Reeb components with complex leaves and their symmetries I : The automorphism groups and Schr\"{o}der's equation on the half line, preprint.\\
{[Hu]} S. Hurder, The $\overline{\partial}$-operator, Appendix A,  In: C. Moore, C.  Schochet, Global analysis on foliated spaces. 
Second edition. Mathematical Sciences Research Institute Publications, {\bf 9}, Cambridge University Press, New York, 2006.\\
{[JS]} W. Jaco, P. Shalen, Seifert fibered spaces in 3-manifolds, Memoirs of the American Mathematical Society,{\bf  21} (220), 1979.\\
{[Jo]} K. Johannson, Homotopy equivalences of 3-manifolds with boundaries, Lecture Notes in Mathematics, {\bf 761}. Springer, Berlin, 1979.\\
{[KH]}  A. Katok, B. Hasselblatt, Introduction to the modern theory of dynamical systems, With a supplementary chapter by Katok and Leonardo Mendoza. Encyclopedia of Mathematics and its Applications, \textbf{54}, Cambridge University Press, Cambridge, 1995. \\
{[Ko1]} F. Kopei,  A remark on a relation between foliations and number theory,  Foliations 2005,  245--249, World Sci. Publ., Hackensack, NJ, 2006.\\
{[Ko2]} F. Kopei, A foliated analogue of one- and two-dimensional Arakelov theory, Abh. Math. Semin. Univ. Hambg.  {\bf 81}  (2011),  no. 2, 141--189. \\
{[Mi]} T. Mihara,  Cohomological approach to class field theory in arithmetic topology, to appear in Canadian J. Math.\\
{[Mo]} M. Morishita, Knots and Primes -- An Introduction to Arithmetic Topology, Universitext, Springer, 2011.\\
{[MT]} M. Morishita, Y. Terashima, Geometry of polysymbols, Math. Res. Lett.  {\bf 15}  (2008),  no. 1, 95--115.\\
{[NN]} A. Newlander, L. Nirenberg, Complex analytic coordinates in almost complex manifolds, 
Ann. of Math. (2)  {\bf 65}  (1957), 391--404. \\
{[NU]} H. Niibo, J. Ueki, Id\`{e}lic class field theory for 3-manifolds and very admissible links, to appear in Transactions of Amer. Math. Soc. \\
{[Sm]} S. Smale, Diffeomorphisms with many periodic points, In: Differential and Combinatorial Topology 1965 (A Symposium in Honor of Marston Morse)  63--80 Princeton Univ. Press, Princeton, N.J. \\
{[St]} J. Stelzig, Weil-reziprozit\"{a}t f\"{u}r Bl\"{a}tterungen durch Riemannsche Fl\"{a}chen, Masterarbeit, Westf\"{a}lische Wilhelms-Universit\"{a}t, M\"{u}nster.\\
{[Ta]} I. Tamura, Topology of foliations: an introduction. Translated from the 1976 Japanese edition and with an afterword by Kiki Hudson. With a foreword by Takashi Tsuboi. Translations of Mathematical Monographs, {\bf 97}. Amer. Math. Soc., Providence, RI, 1992.\\
{[Te]} Y. Terashima, Higher dimensional parallel transports for Deligne cocycles, Integrable systems, topology, and physics (Tokyo, 2000), 291--312, Contemp. Math., {\bf 309}, Amer. Math. Soc., Providence, RI, 2002.
\\
{[Ti]} D. Tischler, On fibering certain foliated manifolds over $S^1$, Topology, {\bf 9} (1970), 153--154.\\
\end{flushleft}
\vspace{.2cm}
{\small
J. Kim:\\
Graduate School of Mathematics, Kyushu University, \\
744, Motooka, Nishi-ku, Fukuoka  819-0395, Japan.\\
e-mail: : j-kim@math.kyushu-u.ac.jp\\
\\
M. Morishita:\\
Graduate School of Mathematics, Kyushu University, \\
744, Motooka, Nishi-ku, Fukuoka  819-0395, Japan.\\
e-mail: morisita@math.kyushu-u.ac.jp \\
\\
T. Noda:\\
Department of Mathematics, Toho University\\
2-2-1, Miyama, Funabashi-shi, Chiba 274-8510, Japan.\\
e-mail: noda@c.sci.toho-u.ac.jp\\
\\
Y. Terashima:\\
Department of Mathematics, Tohoku University,\\
 6-3, Aoba, Aramaki-aza, Aoba-ku, Sendai-shi, Miyagi  980-8578, Japan.\\
e-mail: yujiterashima@tohoku.ac.jp}

\end{document}